\documentclass[11pt]{article}
\usepackage[utf8]{inputenc}
\usepackage{mathrsfs}
\usepackage{mathtools}
\usepackage{enumerate}
\usepackage{amsmath}
\usepackage{amsthm}
\usepackage{amssymb}
\usepackage{multicol}
\usepackage[unicode=true,bookmarks=false,breaklinks=false,pdfborder={0 0 1},backref=section,colorlinks=true]{hyperref}
\hypersetup{citecolor=blue,urlcolor=blue,linkcolor=blue}

\makeatletter
\usepackage{graphicx}
\usepackage{epstopdf}
\usepackage{cite}
\usepackage{todonotes}
\usepackage{subcaption}
\usepackage[font=small]{caption}
\usepackage{xcolor}
\usepackage{fullpage}
\usepackage[misc]{ifsym}
\usepackage{comment}
\usepackage{fontawesome5}
\usepackage{tikz-cd}
\usepackage{algorithm}
\usepackage[noend]{algpseudocode}
\usepackage{lipsum}
\usepackage{multirow}
\usepackage{makecell}
\providecommand{\U}[1]{\protect\rule{.1in}{.1in}}
\newtheorem{theorem}{Theorem}[section]

\newtheorem{lemma}{Lemma}[section]
\newtheorem{definition}{Definition}[section]

\newtheorem{corollary}{Corollary}[section]
\newtheorem{remark}{Remark}[section]
\newtheorem{example}{Example}[section]
\newtheorem{proposition}{Proposition}[section]

\numberwithin{equation}{section}
\numberwithin{equation}{section}
 
\DeclareMathOperator*{\argmin}{\arg\!\min}

\begin{document}
\title{Weak Quasi-Convexity: First- and Second-Order Conditions, Monotonicity, and Quadratic Functions}
\author{Soumitra Dey$^{1,*}$, Simeon Reich$^{2}$\bigskip{}\\
$^{1}$ Einstein Institute of Mathematics, Hebrew University of Jerusalem, 
\\Jerusalem 9190401, Israel \medskip\\
$^{2}$ Department of Mathematics, The Technion -- Israel Institute of Technology, 
\\Haifa 3200003, Israel \medskip\\
}
\date{\today}
\maketitle

\begingroup
\renewcommand{\thefootnote}{\fnsymbol{footnote}}
\footnotetext[0]{%
\begin{minipage}{\textwidth}
\small
\noindent
Soumitra Dey\\
\Letter\ deysoumitra2012@gmail.com

\medskip
\noindent
Simeon Reich\\
\Letter\ sreich@technion.ac.il

\medskip
\noindent
$^{*}$ Corresponding author
\end{minipage}
}
\endgroup


\begin{abstract}
In this paper we introduce weakly quasi-convex functions, which are a generalization of weakly convex functions. 
We study the behavior of weakly quasi-convex functions in Euclidean spaces. 
We provide several illustrative examples that demonstrate the scope and distinctive features of this new class of generalized convex functions. 
We derive a necessary condition for weak quasi-convexity under certain assumptions and provide an example showing that its converse is not true for discontinuous functions. 
Under the continuity assumption on the function involved, a sufficient condition is provided with an adjusted weak quasi-convexity modulus.
Also, we bring out various kinds of monotonicity in connection with weakly quasi-convex functions through the gradient mapping. 
We provide several illustrative examples that interpret and emphasize the newly introduced notions.
We establish a first-order necessary condition for the class of differentiable weakly quasi-convex functions. 
Conversely, we show that a differentiable function satisfying this first-order condition is weakly quasi-convex with an adjusted modulus.
We provide a second-order necessary condition under appropriate assumptions using second-order Dini derivatives. 
As a consequence, we derive a Hessian-based necessary condition for the class of twice continuously differentiable weakly quasi-convex functions.
We provide an example to show that the converse of the second-order characterization does not hold. 
Furthermore, we establish various kinds of necessary conditions concerning general quadratic functions and provide many counterexamples for the failure of their sufficiency. Finally, we derive several consequences of weak quasi-convexity for optimization problems.
\end{abstract}
 
\medskip{}


\section{Introduction}
\label{Sec:1}

\subsection{Background and motivation}
The theory of convexity plays a crucial role in the field of operations research (see, for instance, \cite{HHBA2017}) in finite- and infinite-dimensional Hilbert spaces or, more generally, in normed linear spaces. 
Particularly, in convex optimization including its associated bilevel theory, convex functions constitute one of the fundamental classes of functions. 
A substantial theory has been developed by various researchers (see, for instance, \cite{ABEC2017}) over the years. 
At the same time, various generalizations and related notions have been introduced and extensively studied, including
strongly convex functions, quasi-convex functions, strongly-quasi convex functions, and many others.

Another important generalization of convex functions is the class of weakly convex functions, which play a significant role in the context of non-convex optimization. 
Indeed, a function is said to be weakly convex if the addition of a suitably scaled squared norm makes it convex. 
This transformation property constitutes a useful intermediate framework between convex and non-convex optimization and allows many analytical and algorithmic tools from convex analysis to be extended to a broader class of problems.
Since the theory of non-convex optimization not been developed significantly yet due to its inherent complexity, researchers often assume that the cost function is weakly convex.
In practice, since convexity is too restrictive, considerable attention has been received in developing optimization theory that constitutes weaker forms of convexity.  
In this regard, a new class of convex functions was introduced and studied in \cite{HHUA2025}, and further developments can be found in \cite{MHAl2026}.

Recent studies have also demonstrated the usefulness of generalized convexity in algorithmic optimization. 
For instance, in \cite{FYLI2025}, studied a proximal point algorithm for the optimization problem of a non-smooth, weakly convex function on Euclidean spaces and the order of convergence was computed. Since one of the fundamental challenges in computer science is to solve complex problems efficiently in terms of both time and space complexity, computing the order of convergence attracted considerable attention from researchers. 
The class of strongly quasi-convex functions was introduced in \cite{BTPO1966}, and a substantial theory has been developed for the class of strongly quasi-convex functions in \cite{SMGR2025,NMNA2025}. 

Motivated by the above developments over the years, in this paper we have introduced a new class of functions, called ``\textit{weakly quasi-convex functions}", which provides a flexible framework for further extending the theory of generalized convexity. 
This class is sufficiently broad to include several important classes of generalized convex functions while retaining structural properties that are useful for optimization. 
The relationships between weakly quasi-convex functions and other relevant classes of convex functions are illustrated in Figure \ref{fig:convexity-relations}.

A notion of weak quasi-convexity is presented in \cite[Definition 2.1]{MHAR2018}.
Specifically, it says that an objective function $f$ is $\eta$-weakly quasi-convex on a domain $C$ with respect to a global minimizer $x^*$ if there exists a constant $\eta>0$ such that, at every point $x\in C$, the gradient $\nabla f(x)$ has a sufficiently large component in the direction $x-x^*$, bounded below by a positive constant $\eta$ times the objective-value gap $f(x)-f(x^*)$. 
The notion of weak quasi-convexity we have introduced in this paper is defined through a quadratic relaxation of the classical quasi-convexity inequality, which is different from the one studied in \cite[Definition 2.1]{MHAR2018}.

\subsection{Our contributions}
The main contributions of this paper are presented below.

\begin{enumerate}
\item We have introduced a new class of generalized convex functions, which we call weakly quasi-convex functions by incorporating a quadratic error term into the usual quasi-convexity inequality on Euclidean spaces. 
The relationships between this class and several classical generalized convexity classes are investigated.

\item We have derived a new necessary condition for weak quasi-convexity of a function under appropriate assumptions on Euclidean space. 
We have shown by counterexample that the converse fails for discontinuous functions and that the segmentwise ordering assumption cannot, in general, be replaced with an appropriate direction. 
However, a converse result yielding weak quasi-convexity with an adjusted modulus is established under the continuity assumption of the function involved in the analysis.

\item We have constructed several nontrivial examples and counterexamples to distinguish weak quasi-convexity of a function from various kinds of existing convexity notions such as weak convexity, quasi-convexity, and strong convexity. 
Particularly, an example shows that the new class contains functions that are neither weakly convex nor quasi-convex in the $1$-dimensional Euclidean space.

\item We have developed a first-order characterization for the class of differentiable weakly quasi-convex functions and an example is provided to illustrate the theoretical result. 
A necessary gradient inequality is established, and together with a converse result showing that the corresponding first-order condition implies weak quasi-convexity with an adjusted modulus.

\item We have introduced new generalized monotonicity notions, namely weak quasi-monotonicity, weak pseudo-monotonicity, and strict weak pseudo-monotonicity. 
Their relationships with classical quasi-monotonicity and pseudo-monotonicity are analyzed through examples and counterexamples, which distinguish our notions of monotonicity with the relevant existing one.

\item For the class of differentiable functions, a connection is drawn between weak quasi-convexity and weak quasi-monotonicity through the gradient mapping. 
In particular, we have shown that the gradient of a differentiable weakly quasi-convex function is weakly quasi-monotone. We have provided a counterexample for the converse of this result.

\item We have established a second-order characterization for the weakly quasi-convex function under appropriate assumptions using the second-order Dini derivative. 
From this result, we have derived a corollary as a second-order necessary characterization for the class of twice continuously differentiable weakly quasi-convex functions. 
Specifically, a Hessian condition is obtained along directions orthogonal to the gradient, and its limitations are demonstrated through a counterexample.

\item We have studied several stability and closure-type properties of weakly quasi-convex functions, including scalar multiplication, composition with linear mappings, pointwise maxima, composition with suitable nondecreasing functions, and one-dimensional restrictions along line segments.

\item We have analyzed the weak quasi-convexity of quadratic functions in detail. 
Necessary and sufficient conditions have been obtained on the whole Euclidean space, and relationships among weak convexity, weak quasi-convexity, and spectral properties of the associated symmetric matrix have been established.
Additionally, quadratic functions on restricted domains are also investigated. 
We have obtained necessary conditions involving the number of negative eigenvalues of the matrix \(Q+\rho I\), together with counterexamples showing the failure of the corresponding converses.

\item An optimization consequence of weak quasi-convexity is obtained. A necessary first-order inequality satisfied by a global minimizer over a convex feasible set is derived, and its connection with a variational inequality problem is discussed.
An application of the first-order necessary condition to the saddle point problem is provided.
Additionally, we have established a result on the necessary and sufficient characterization of weak quasi-convexity by drawing a connection with the approximate solution set of a constrained optimization problem under the assumption that the associated function $f$ is bounded below.
\end{enumerate}


\subsection{Notations}
Denote
$\mathbb{R}^n_+=\{x\in\mathbb{R}^n\;\vert\; x_i\geq 0\;\forall 1\leq i\leq n\}$, $\mathbb{R}^n_{++}=\{x\in\mathbb{R}^n\;\vert\; x_i> 0\;\forall 1\leq i\leq n\}$, and $e_i$ to represent a vector in $\mathbb{R}^n$ with $i$-th coordinate is $1$. 
The notations $\langle x,y\rangle=x^{\top}y$ and $\|x\|=\sqrt{\langle x,x\rangle}$ denote the standard Euclidean inner product and norm on $\mathbb{R}^n$, respectively.
For a nonempty subset $C$ of $\mathbb{R}^n$, int$(C)$ denotes the interior of $C$ and for $x\in\mathbb{R}^n$, $B(x,\delta(x))$ denotes the open ball centered at $x$ with radius $\delta(x)$.
The level set of $f:\mathbb{R}^n\rightarrow\mathbb{R}$ at height $\alpha\in\mathbb{R}$ is 
$L(f,\alpha)=\{x\in\mathbb{R}^n\vert\; f(x)\leq \alpha\}$ (see, \cite[p. 6]{HHBA2017}). 
The notation $\subset$ stands for set inclusion, while $\subsetneq$ stands for proper set inclusion.
For a real-valued differentiable function $f$ on an open subset of $\mathbb{R}^n$, $\nabla f(x)$ denotes its gradient at $x$ and $\nabla^2 f(x)$ denotes its Hessian matrix at $x$.
A matrix $M\in\mathbb{R}^{n\times n}$ is symmetric if $M=M^T$, transpose of $M$. A matrix $M$ is said to be positive semidefinite (psd) if it is symmetric and $x^TMx\geq 0$, for all $x\in\mathbb{R}^n$. Let $I$ be the identity matrix and $M\geq 0$ denote that all the entries in $M$ are nonnegative. 
For a matrix $M<0$ means that all elements of it are negative. The notation $(x,y]=\{z\;\vert\; z=\lambda y+(1-\lambda)x,\;\lambda\in (0,1]\}$. 
For an interval $I\subset \mathbb{R}$, denote
\begin{eqnarray*}
    L^{1}(I):=\left\{f:I\to\mathbb{R}\;\vert\;f \text{ is Lebesgue measurable and }\int_I |f(x)|\,dx<\infty\right\}.
\end{eqnarray*}

\subsection{The organization of this paper} 
This paper is organized as follows.
In Section \ref{Sec:2}, we recall some basic definitions and preliminary results. 
In Section \ref{Sec:3}, we introduce the notions of weak quasi-convexity and strict weak quasi-convexity, together with other relevant concepts.
We derive a first-order necessary condition for the class of weakly quasi-convex functions and establish its converse with an appropriate adjustment of the weak quasi-convexity modulus. 
We also introduce and study several generalized monotonicity concepts.
Furthermore, we obtain a second-order necessary condition through the second-order Dini derivative, together with several relevant consequences.
In Section \ref{Sec:4}, we derive further properties of weakly quasi-convex functions and analyze weak quasi-convexity of quadratic functions.
In Section \ref{Sec:OptimizationProblem}, we establish a necessary optimality condition for optimization problems involving weakly quasi-convex functions and discuss its connection with a variational inequality problem.
Additionally, we give an application to the saddle point problem.
Also, we derive a connection of approximate solution sets of a constrained optimization problem with the weak quasi-convexity of the associated function involved in the constrained optimization.
Finally, in Section \ref{Sec:6}, we draw a conclusion of our work.

\section{Preliminaries}
\label{Sec:2}

In this section, we recall several preliminary definitions and results that we will be using in the sequel.

\begin{definition}\label{Def1:convex}\cite[p. 137]{ABEC2023}
    Let $C$ be a nonempty convex subset of $\mathbb{R}^n$. A function $f:C\rightarrow\mathbb{R}$ is said to be
    \begin{enumerate}
        \item convex on $C$ if for all $x,y\in C$ and $\lambda\in [0,1]$ such that
        \begin{equation}\label{ineq1:Def1:convex} 
        f(\lambda y+(1-\lambda)x)\leq\lambda f(y)+(1-\lambda)f(x).
        \end{equation}
        
        \item strictly convex on $C$ if for all $x,y\in C$ with $x\neq y$ and $\lambda\in (0,1)$ such that
        \begin{equation}\label{ineq1:Def1:convex} 
        f(\lambda y+(1-\lambda)x)<\lambda f(y)+(1-\lambda)f(x).
        \end{equation}        
    \end{enumerate}
\end{definition}

\begin{definition}\label{Def3:quasi-convex}(\cite[p. 161]{ABEC2023}\cite[p. 184]{HHBA2017})
    Let $C$ be a nonempty convex subset of $\mathbb{R}^n$. A function $f:C\rightarrow\mathbb{R}$ is said to be
    \begin{enumerate}
        \item quasi-convex on $C$ if for all $x,y\in C$ and $\lambda\in [0,1]$ such that
        \begin{equation}\label{ineq1:Def1:convex} 
        f(\lambda y+(1-\lambda)x)\leq\max\{f(y),f(x)\}.
        \end{equation}
        
        \item strictly quasi-convex on $C$ if for all $x,y\in C$ with $x\neq y$ and $\lambda\in (0,1)$ such that
        \begin{equation}\label{ineq1:Def1:convex} 
        f(\lambda y+(1-\lambda)x)<\max\{f(y),f(x)\}.
        \end{equation}        
    \end{enumerate}
\end{definition}

\begin{definition}\label{Def3:Strongly_quasi-convex}
     Let $C$ be a nonempty convex subset of $\mathbb{R}^n$. A function $f:C\rightarrow\mathbb{R}$ is said to be
     \begin{enumerate}
         \item strongly convex with modulus $\rho>0$ if for all $x,y\in C$ and $\lambda\in [0,1]$ such that
            \begin{equation}\label{ineq7} 
              f(\lambda y+(1-\lambda)x)\leq\lambda f(y)+(1-\lambda)f(x)-\frac{\rho}{2}\lambda(1-\lambda)\|y-x\|^2.
             \end{equation}
        \item strongly quasi-convex with modulus $\rho>0$ if for all $x,y\in C$ and $\lambda\in [0,1]$ such that
            \begin{equation}\label{ineq7} 
              f(\lambda y+(1-\lambda)x)\leq\max\{f(y), f(x)\}-\frac{\rho}{2}\lambda(1-\lambda)\|y-x\|^2.
             \end{equation}     
     \end{enumerate}
\end{definition}

\begin{remark}\label{rmk1:Def3:Strongly_quasi-convex}
    The notion of Definition \ref{Def3:Strongly_quasi-convex} is borrowed from \cite[Definition 4.1]{JPVI1983}.
\end{remark}

\begin{definition}\cite[Definition 4.1]{JPVI1983}
Let $C$ be a nonempty convex subset of $\mathbb{R}^n$. A function $f:C\rightarrow\mathbb{R}$ is said to be $\rho$-weakly convex if there exists a constant $\rho>0$ such that for all $x,y\in C$ and $\lambda\in [0,1]$
\begin{equation}\label{ineq7} 
f(\lambda y+(1-\lambda)x)\leq\lambda f(y)+(1-\lambda)f(x)+\frac{\rho}{2}\lambda(1-\lambda)\|y-x\|^2.
\end{equation}
\end{definition}

\begin{remark}\label{rmk1:weak:conv:nessuff}
 It is also well known that the function $f$ defined on a nonempty convex subset $C$ of $\mathbb{R}^n$ is $\rho$-weakly convex on $C$ if and only if $f(\cdot)+\frac{\rho}{2}\|\cdot\|^2$ is convex on $C$. For example, one can see \cite{FYLI2025}, for an extended real-valued function on $\mathbb{R}^n$.
\end{remark}

\begin{lemma}\label{Lem:Quasi-convex}\cite[Definition 7.33]{ABEC2023}
 A function $f:\mathbb{R}^n \rightarrow \mathbb{R}$ is quasi-convex if and only if all level sets of $f$ are convex.
\end{lemma}

\begin{lemma}\cite[Theorem 7.36]{ABEC2023}\label{Lemma:conv&continuity}
    Let $C$ be a nonempty convex subset of $\mathbb{R}^n$. A convex function $f:C\rightarrow\mathbb{R}$ is continuous on int$(C)$.
\end{lemma}

\begin{lemma}\label{lem1:conv&continuity}
   Let $C$ be a nonempty convex subset of $\mathbb{R}^n$. If $f:C\rightarrow\mathbb{R}$ is weakly convex with modulus $\rho>0$ on $C$, then it is continuous on int$(C)$. 
\end{lemma}

\begin{proof}
    Let $g(x):=f(x)+(\rho/2)\|x\|^2$ for $x\in C$. By the definition of weak convexity of $f$, the function $g$ is convex on $C$. Using Lemma \ref{Lemma:conv&continuity}, we get $g$ is continuous on $\text{int}(C)$. Therefore, $f(x)=g(x)-(\rho/2)\|x\|^2$ is continuous on $\text{int}(C)$.
\end{proof}

Let $C$ be a nonempty, open and convex subset of $\mathbb{R}^n$ and $f:C\rightarrow\mathbb{R}$ be a function. We say that the function $f$ has a directional derivative at $x\in C$ in the direction $d\in\mathbb{R}^n$ if 
\begin{eqnarray}\label{ineq1:Directional_Deri}
    f^\prime(x;d)=\lim_{h\to 0^+}\frac{f(x+hd)-f(x)}{h}
\end{eqnarray}
exists.

We say that the function $f$ has the second-order Dini derivative on $C$ if
\begin{eqnarray}\label{ineq1:Dini_second}
    D^{\prime\prime} f(x;d)=\lim_{h\to 0^+}\frac{f^\prime(x+hd;d)-f^\prime(x;d)}{h}
\end{eqnarray}
exists. We refer the reader to \cite{SSHI1993, RCOM1990, WLCH1994} for more details on this topic.

\begin{definition}
    Let $I$ be an interval of $\mathbb{R}$. A function $g:I\rightarrow\mathbb{R}$ is said to be absolutely continuous on $I$ if for any $\epsilon>0$ there exists a $\delta>0$ such that whenever a finite sequence of disjoint open intervals $\{(x_k,y_k)\}_{k=1}^\ell$ satisfies 
    \begin{eqnarray}
       \sum_{k=1}^\ell (y_k-x_k)<\delta\implies \sum_{k=1}^\ell \vert g(y_k)-g(x_k)\vert<\epsilon.
    \end{eqnarray}
\end{definition}

\noindent It is well known that if $g$ is an absolutely continuous function on a compact interval $I\subset\mathbb{R}$, then 
\begin{enumerate}
    \item[(i)] it is almost everywhere (a.e.) differentiable on $I$.
    \item[(ii)] $g^\prime\in L^1(I)$.
    \item[(iii)] $g(s)-g(t)=\int_{t}^s g^\prime (\omega)\,d\omega$, for all $s,t\in I$ with $t\leq s$. 
\end{enumerate}






\section{Weak quasi-convexity}
\label{Sec:3}

\subsection{Weak quasi-convexity and strict weak quasi-convexity}

\begin{definition}\label{DEF:MAIN;WQC}
Let $C$ be a nonempty convex subset of $\mathbb{R}^n$. A function $f:C\rightarrow \mathbb{R}$ is said to be 
\begin{enumerate}
 \item[(1)] $\rho$-\textit{weakly quasi-convex} or \textit{weakly quasi-convex} on $C$ if there exists a $\rho>0$ such that for all $x,y\in C$ and $\lambda\in [0,1]$
\begin{equation}\label{our:Def:1}
 \begin{aligned}
  f(\lambda y+(1-\lambda)x)\leq\max\{f(y),f(x)\} + \frac{\rho}{2}\lambda(1-\lambda)\|y-x\|^2.
 \end{aligned}
\end{equation}

\item[(2)] $\rho$-\textit{strictly weakly quasi-convex} or \textit{strictly weakly quasi-convex} on $C$ if there exists a $\rho>0$ such that for all $x,y\in C$ with $x\neq y$ and $\lambda\in (0,1)$
\begin{equation}\label{ineq1:DEF:MAIN;WQC}
\begin{aligned}
f(\lambda y+(1-\lambda)x)<\max\{f(y),f(x)\}+\frac{\rho}{2}\lambda(1-\lambda)\|y-x\|^2.
\end{aligned}
\end{equation}
\end{enumerate}
\end{definition}

\begin{theorem}\label{thm1:CharWQC}
        Let $C$ be a nonempty convex subset of $\mathbb{R}^n$. Then the function $f:C\rightarrow\mathbb{R}$ is $\rho$-weakly quasi-convex on $C$ if and only if for any $x,y\in C$ and $z\in [x,y]\subset C$ such that
    \begin{equation}\label{ineq1:thm1:CharWQC}
        \begin{aligned}
           f(y)\leq f(x) ~ \Longrightarrow ~ f(z) \leq f(x)+\frac{\rho}{2}\|z-y\|\|z-x\|.
        \end{aligned}
    \end{equation}
\end{theorem}

\begin{proof}
    Since $z\in [x,y]$, there exists $\lambda\in [0,1]$ such that $z=\lambda y+(1-\lambda)x$. The result follows using the fact that $\|z-y\|=(1-\lambda)\|y-x\|$ and $\|z-x\|=\lambda\|y-x\|$.
\end{proof}

A similar result holds for $\rho$-strictly weakly quasi-convex functions presented below.

\begin{theorem}\label{thm1:CharSWQC}
        Let $C$ be a nonempty convex subset of $\mathbb{R}^n$. Then the function $f:C\rightarrow\mathbb{R}$ is $\rho$- strictly weakly quasi-convex on $C$ if and only if for any $x,y\in C$ with $x\neq y$ and $z\in (x,y)\subset C$ such that
    \begin{equation}\label{ineq1:thm1:CharSWQC}
        \begin{aligned}
           f(x)\geq f(y) ~ \Longrightarrow ~ f(x) > f(z)-\frac{\rho}{2}\|z-y\|\|z-x\|.
        \end{aligned}
    \end{equation}
\end{theorem}

\begin{proof}
    The proof is immediate from the proof of Theorem \ref{thm1:CharWQC}.
\end{proof}



\begin{remark}
 In particular, when $\rho=0$, Definition \ref{DEF:MAIN;WQC} reduces to the classical notions of quasi-convex and strictly quasi-convex functions on $C$, respectively.
\end{remark}

\begin{figure}[htbp]
\centering
\begin{tikzpicture}[
    node distance=0.3cm and 0.3cm,
    box/.style={
        draw,
        rectangle,
        rounded corners,
        minimum width=0.5cm,
        minimum height=0.1cm,
        align=center,
        font=\small
    },
    arrow/.style={-{Latex[length=1.3mm]}}
]
\node[box] (StronglyConvex) at (-5,0) {Strongly\\ convex};

\node[box] (StrictlyConvex) at (-2,2) {Strictly\\ convex};
\node[box] (Convex) at (2,2) {Convex};
\node[box] (WeaklyConvex) at (6,2) {Weakly\\ convex};

\node[box] (StrictlyQuasiConvex) at (-2,0) {Strictly\\ quasi-convex};
\node[box] (QuasiConvex) at (2,0) {Quasi-convex};
\node[box] (StrictlyWeaklyQuasiConvex) at (6,0) {Strictly weakly\\ quasi-convex};

\node[box] (StronglyQuasiConvex) at (-2,-2) {Strongly\\ quasi-convex};
\node[box] (WeaklyQuasiConvex) at (6,-2) {Weakly\\ quasi-convex};

\draw[arrow] (StronglyConvex) -- (StrictlyConvex);
\draw[arrow] (StronglyConvex) -- (StronglyQuasiConvex);

\draw[arrow] (StrictlyConvex) -- (Convex);
\draw[arrow] (Convex) -- (WeaklyConvex);

\draw[arrow] (StrictlyConvex) -- (StrictlyQuasiConvex);
\draw[arrow] (StronglyQuasiConvex) -- (StrictlyQuasiConvex);

\draw[arrow] (StrictlyQuasiConvex) -- (QuasiConvex);
\draw[arrow] (QuasiConvex) -- (StrictlyWeaklyQuasiConvex);

\draw[arrow] (WeaklyConvex) -- (StrictlyWeaklyQuasiConvex);
\draw[arrow] (StrictlyWeaklyQuasiConvex) -- (WeaklyQuasiConvex);

\draw[arrow] (Convex) -- (QuasiConvex);
\draw[arrow] (Convex) -- (StrictlyQuasiConvex);

\end{tikzpicture}
\caption{Relations among different classes of convex functions.}
\label{fig:convexity-relations}
\end{figure}
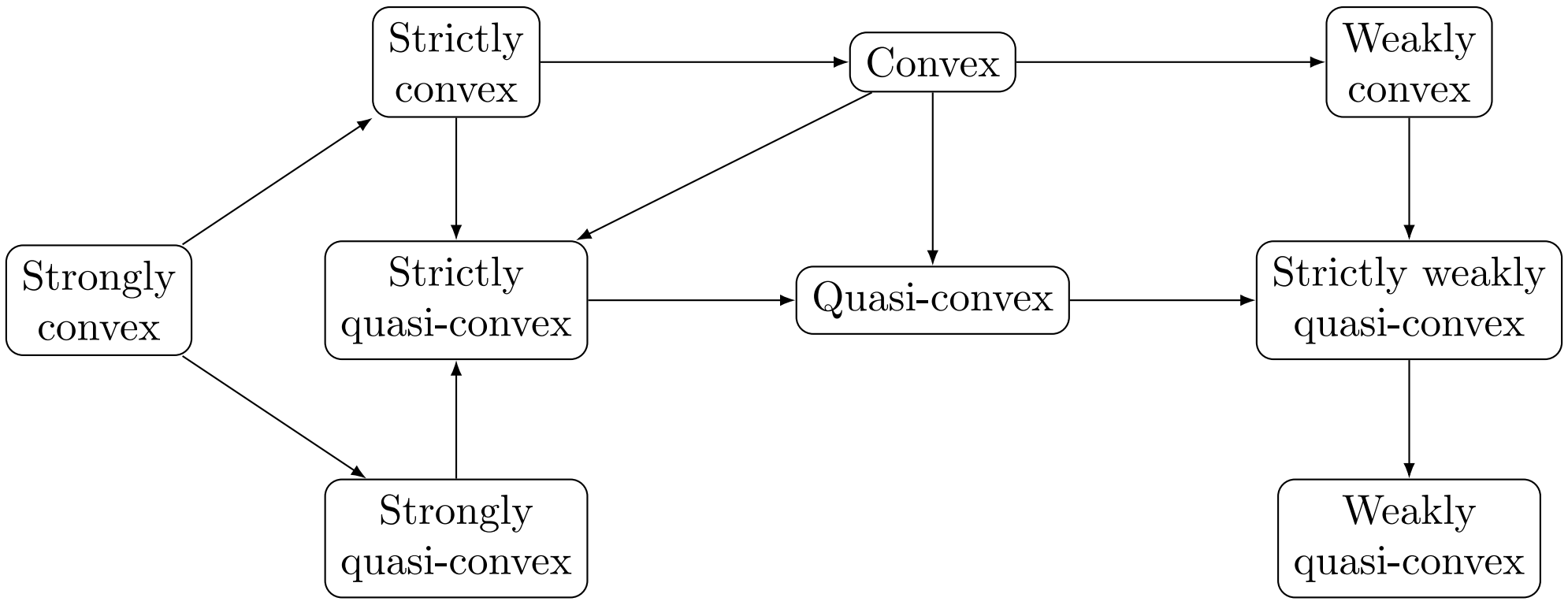

\begin{remark}\label{rmk1}
 Clearly, every weakly convex function is weakly quasi-convex. Furthermore, Figure \ref{fig:convexity-relations} shows that the class of weakly quasi-convex functions is larger than several other classes of convex functions on a nonempty convex subset $C$ of $\mathbb{R}^n$.
Observe that the two classes associated to quasi-convex and weakly convex functions are independent. 
\end{remark}

The following example shows that quasi-convexity does not guarantee the weak convexity of a function.

\begin{example}\label{example:1}
 Let us consider the function  $f: \mathbb{R} \rightarrow \mathbb{R}$ defined by
\[
f(x)=
\begin{cases}
1,   & x<0\\
x^2, & x\geq 0.
\end{cases}
\]

The level sets of $f$, $L(f,\alpha)=\{x\in\mathbb{R}\;\vert\;f(x)\leq\alpha\},\;\alpha\in\mathbb{R},$
are convex.  Specifically, for $\alpha<0$, $L(f,\alpha)=\varnothing,$ convex. 
For $\alpha=0$, $L(f,\alpha)=\{0\}$, convex.
For $0<\alpha< 1$, if $x\geq 0$, $L(f,\alpha)=[0,\sqrt{\alpha}],$ 
convex.
For each $\alpha\geq 1$, if $x<0$, $L(f,\alpha)=(-\infty,0)$, and if $x\geq 0$, $L(f,\alpha)=[0,\sqrt{\alpha}]$. Therefore,
for $\alpha\geq 1$, $L(f,\alpha)=(-\infty,\sqrt{\alpha}],$ convex.
Therefore,  $f$ is a quasi-convex function.

Clearly, $f$ is discontinuous at $0.$ Therefore, by Lemma \ref{Lemma:conv&continuity}, $f$ is not convex. Moreover, from Lemma \ref{lem1:conv&continuity}, $f$ is not weakly convex on $\mathbb{R}$.
\end{example}

The following example shows that weak convexity does not guarantee quasi-convexity of a function.

\begin{example}\label{example:2}
 Let us consider a function $f: \mathbb{R} \rightarrow \mathbb{R}$ by $f(x) = x^2 + \sin(5x)$. Consider the function $g(x) := (1 + \rho/2)x^2 + \sin(5x)$, with $\rho > 0$. Clearly, the second derivative of $g$ is non-negative for any $\rho\geq 23$. 
Therefore, the function $f$ is weakly convex for any $\rho\geq 23$. 
However, it is not quasi-convex since its level set at height $\alpha=1$ is not convex, see the geometry given in Figure \ref{fig:quasi_convex_level_set}.

\begin{figure}[htbp]
    \centering
    \includegraphics[width=0.35\textwidth]{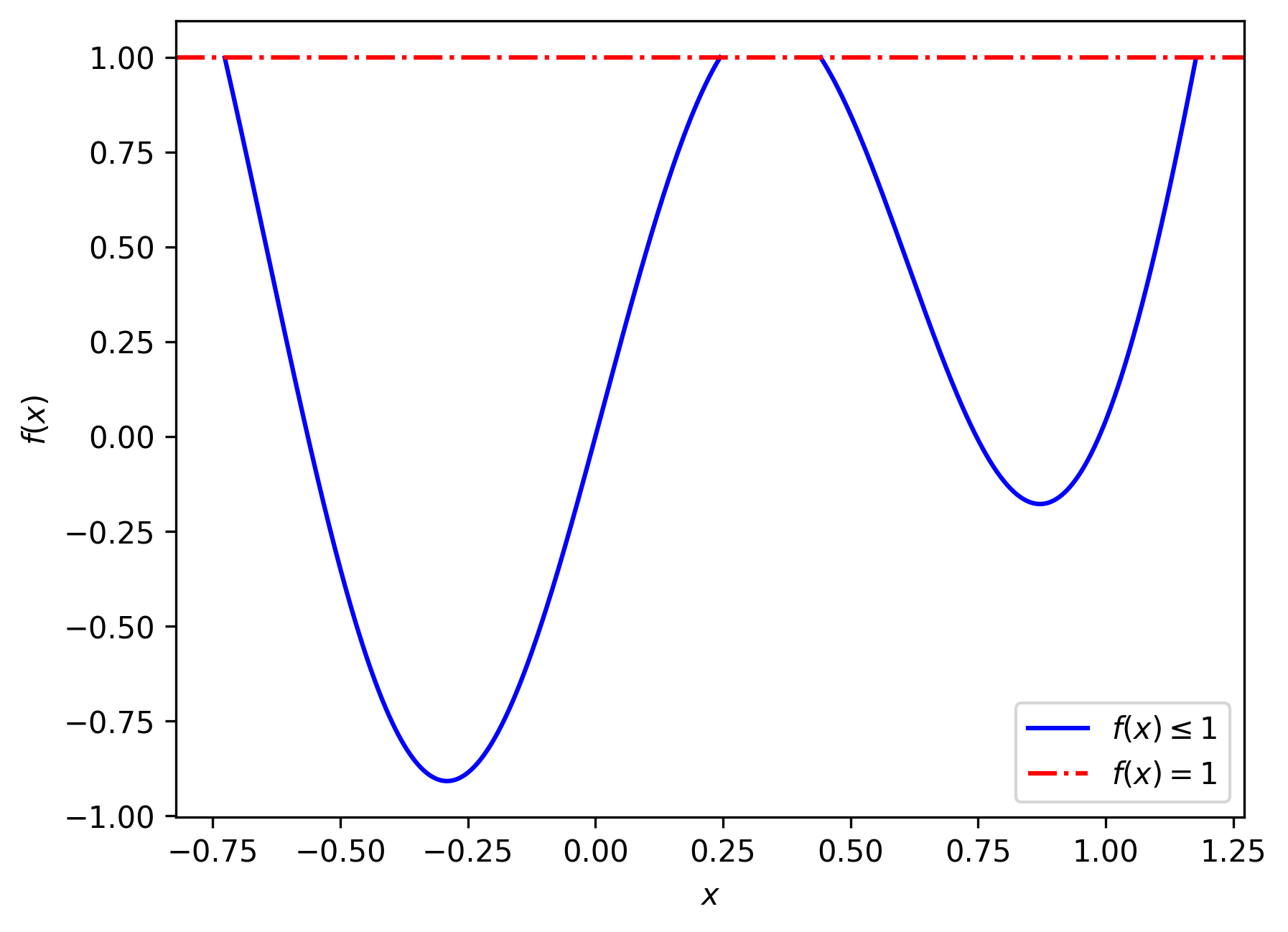}
    \caption{$L(f,1)=\{x\in\mathbb{R}:x^2+\sin(5x)\leq 1\}$}
    \label{fig:quasi_convex_level_set}
\end{figure}
\end{example}

Next we present a lemma which will be useful for our next example.

\begin{lemma}\label{Lem1:NextExample}
Let $g:[0,1]\rightarrow\mathbb{R}$ be a continuous function and twice continuously differentiable on $(0,1)$. Suppose, there exists a constant $L\geq 0$ such that $g^{\prime\prime}(t)\geq -L$, for all $t\in (0,1)$. Then 
\begin{eqnarray}
    g(t)\leq t g(1)+(1-t)g(0)+\frac{t(1-t)}{2}L,\;\text{for all}\; t\in [0,1].
\end{eqnarray}
\end{lemma}

\begin{proof}
    Let $h(t)=g(t)+(L t^2/2)$ for $t\in [0,1]$. Since $g$ is twice continuously differentiable in $(0,1)$, $h$ is so. 
    Again since $g^{\prime\prime}(t)\geq -L$ for all $t\in (0,1)$, $h^{\prime\prime}(t)=g^{\prime\prime}(t)+L\geq 0$ for all $t\in (0,1)$. Therefore, $h$ is convex on $[0,1]$ using its continuity. Therefore, $h(t)\leq t h(1)+(1-t)h(0)$ for all $t\in [0,1]$, which is equivalent to 
    \begin{eqnarray}
    g(t)\leq t g(1)+(1-t)g(0)+\frac{t(1-t)}{2}L,\;\text{for all}\; t\in [0,1].
\end{eqnarray}
\end{proof}

The following example illustrates that the class of weakly quasi-convex functions is a true generalization of many other classes of convex functions. In particular, in the example given below, we consider a function which is weakly quasi-convex but neither weakly convex nor quasi-convex on the one-dimensional Euclidean space. 

\begin{example}\label{Example:WQCnotQCnotWC}
    Consider a function $f:\mathbb{R}\rightarrow\mathbb{R}$ by
\[
f(x)=
\begin{cases}
    \sin x,   & x \geq 0\\
    1, & x < 0.
\end{cases}
\]
Observe that $f$ is not quasi-convex. Indeed, let $x=7\pi/6$, $y=2\pi+(7\pi/6)$, and $\lambda=1/2.$ Then $f(x)=-1/2$, $f(y)=-1/2$, and $f((x+y)/2)=f(2\pi+(\pi/6))=1/2$. This implies that $f((x+y)/2)>\max\{f(x),f(y)\}$. Hence, $f$ is not quasi-convex. Furthermore, $f$ is not a weakly convex function since it is discontinuous at $0$.

On the other hand, $f$ is $1$-weakly quasi-convex on $\mathbb{R}$. Indeed, let $x,y \in \mathbb{R}$ and $\lambda\in [0,1]$. The following cases may arise:

\medskip{}
    \textbf{Case-I:} If $x,y\geq 0$, then $z_\lambda := \lambda y+(1-\lambda)x\geq 0$, for all $\lambda \in [0,1]$. This implies that $f(z_\lambda) = \sin(\lambda y+(1-\lambda)x)$.

    Let $g(\lambda)=f(z_\lambda)$ for $\lambda \in [0,1]$. Using the fact that $\sin\omega \leq 1$, for all $\omega \in \mathbb{R}$, we get $g^{\prime\prime}(\lambda) \geq -(x-y)^2$, for all $\lambda\in [0,1]$. Therefore, using Lemma \ref{Lem1:NextExample}, we obtain
    \begin{eqnarray}\label{ineq2:Example:WQCnotQCnotWC}
        \begin{aligned}
            \sin(\lambda y+(1-\lambda)x)&\leq \lambda \sin(y) + (1-\lambda) \sin (x) + \lambda(1-\lambda)\frac{(x-y)^2}{2} \\&
            \leq \max\{f(y), f(x)\} +\lambda(1-\lambda) \frac{(x-y)^2}{2}.
        \end{aligned}
    \end{eqnarray}
    Thus, $f$ is $1$-weakly quasi-convex.

\medskip{}
    \textbf{Case-II:} If $x,y<0$, then $z_\lambda := \lambda y+(1-\lambda)x <0$, for all $\lambda \in [0,1].$
    This implies that $f(z_\lambda) = 1 \leq \max\{f(x), f(y)\}$ and therefore, $f$ is $\rho$-weakly quasi-convex for any $\rho > 0$. In particular, take $\rho=1$.

\medskip{}
\textbf{Case-III:} If $x< 0\leq y$, then $f(x)=1$ and $f(y)=\sin y$. Now, we have the following cases:
    \begin{itemize}
     \item[$(i)$] If $z_\lambda:=\lambda y+(1-\lambda)x<0$, for some $\lambda\in [0,1]$, then $f(z_\lambda)=1$. Therefore, 
    \begin{eqnarray}\label{ineq3:Example:WQCnotQCnotWC}
    \begin{aligned}
        f(z_\lambda)\leq\max\{f(y), f(x)\}.
    \end{aligned}
    \end{eqnarray}
    Therefore, $f$ is $\rho$-weakly quasi-convex for any $\rho>0.$ Particularly, take $\rho=1$.

\medskip{}
    \item[$(ii)$] If $z_\lambda:=\lambda y+(1-\lambda)x\geq 0$ for some $\lambda\in [0,1]$, then $f(z_\lambda)\leq 1\leq\max\{f(y), f(x)\}.$ Therefore, $f$ is $\rho$-weakly quasi-convex for any $\rho>0.$ Particularly, take $\rho=1$.
\end{itemize}

\medskip{}
    \textbf{Case-IV:} If $y< 0\leq x$, then weak quasi-convexity follows from \textbf{Case-III}.
\end{example}

\begin{remark}
\begin{enumerate}\label{rmk:Disc_bdd_WQC_fun}
\item In general, it is well known that every real valued convex function on $\mathbb{R}$ is continuous. However, it is not true for the weakly quasi-convex functions. See, for instance, Example \ref{Example:WQCnotQCnotWC}.

\item  It is well known that if a function is bounded above on an unbounded convex subset $C$ with nonempty interior of $\mathbb{R}^n$, then it is not a strongly quasi-convex (see, \cite[p. 155]{MVJO1993}). However, this result is not true for the weakly quasi-convex functions as Example \ref{Example:WQCnotQCnotWC} shows. 
\end{enumerate}
\end{remark}

A function $f:\mathbb{R}^n\rightarrow\mathbb{R}$ is $\rho$-weakly convex if and only if $f(\cdot)+(\rho/2)\|\cdot\|^2$ is convex (see, Remark \ref{rmk1:weak:conv:nessuff}). However, an analogous result is not true for the weakly quasi-convex functions. In other words, it is not true that $f$ is $\rho$-weakly quasi-convex if and only if $f(\cdot)+(\rho/2)\|\cdot\|^2$ is quasi-convex. We give an example below for the illustration.

\begin{example}\label{Example1:WQCnotQCnotWC}
Consider a function $f:\mathbb{R}\rightarrow\mathbb{R}$ by
\[
f(x)=
\begin{cases}
 -\frac{x}{2},   & x \leq 1\\
  -\frac{1}{2}-2(x-1), & x > 1.
\end{cases}
\] 
Since the function $f$ is monotonically decreasing, it is quasi-convex on $\mathbb{R}$. Therefore, it is a $\rho$-weakly quasi-convex function on $\mathbb{R}$. Let $\rho=1$, $x=0.9$, and $y=1.1$. Then $f(x)+(1/2)x^2=-0.045$ and $f(y)+(1/2)y^2=-0.095$. Therefore, $f((x+y)/2)+(1/2)((x+y)/2)^2=0>\max\{-0.045, -0.095\}$. Therefore, $f(x)+(\rho/2)x^2$ is not quasi-convex on $\mathbb{R}.$

Now, let us consider the following function $f:\mathbb{R}\rightarrow\mathbb{R}$ defined by
\[
f(x)=
\begin{cases}
    -\frac{x^2}{2}, & x<\frac{1}{2}\\
    \frac{1}{2}-\frac{x^2}{2}, & x\geq\frac{1}{2}.
\end{cases}
\] 

Then 
\[
f(x)+\frac{x^2}{2}=
\begin{cases}
    0, & x<\frac{1}{2}\\
    \frac{1}{2}, & x\geq\frac{1}{2}.
\end{cases}
\] 
is monotonically increasing on $\mathbb{R}$. Therefore, it is a quasi-convex function on $\mathbb{R}.$ However, the function $f$ is not $1-$weakly quasi-convex function on $\mathbb{R}.$ Indeed, let $\lambda=1/2$, $\rho=1$, $x=0$, and $y=1$. Then $f((x+y)/2)-(\rho/2)\lambda(1-\lambda)(x-y)^2=f(1/2)-(1/8)=1/4>0=\max\{f(x), f(y)\}.$
\end{example}

The result below provides another necessary analogue of Theorem \ref{thm1:CharWQC}. 

\begin{theorem}\label{thm1:CharWQC:Analo:prev}
    Let $C$ be a nonempty convex subset of $\mathbb{R}^n$ and the function $f:C\rightarrow\mathbb{R}$ be a $\rho$-weakly quasi-convex on $C$. Let $x,y\in C$ and $z=\lambda y+(1-\lambda)x$ for $\lambda\in (0,1]$. Assume that $C^{x,y}_z:=\{\omega\in(z,y]\;\vert\; f(x)\leq f(\omega)\}=(z,y]$. Then 
    \begin{equation}\label{ineq1:thm1:CharWQC:Analo:prev}
        \begin{aligned}
          f(z) \leq f(y)+\frac{\rho}{4}(1-\lambda^2)\|y-x\|^2.
        \end{aligned}
    \end{equation}
\end{theorem}

\begin{proof}
    Let $z=\lambda y +(1-\lambda)x$ for $x,y\in C$ and $\lambda\in (0,1]$. For $i=0,1,\dots, n$, consider the points 
    \begin{eqnarray}\label{ineq2:thm1:CharWQC:Analo:prev}
        \widetilde{z}_i=z+\frac{i}{n}(y-z).
    \end{eqnarray}
    By our hypothesis $f(x)\leq f(\widetilde{z}_{i+1})$ for all $i=0,1,2,\dots, n-1$. For $i=0,1,2,\dots, n-1$, denote
    \begin{eqnarray}\label{ineq3:thm1:CharWQC:Analo:prev}
        \widetilde{z}_{i+1}=z+\frac{i+1}{n}(y-z)=x+(\lambda+\frac{i+1}{n}(1-\lambda))(y-x):=x+\alpha_{i+1}(y-x),
    \end{eqnarray}
    where $\alpha_{i+1}=\lambda+\frac{i+1}{n}(1-\lambda).$ Then it is easy to check that 
    \begin{eqnarray}\label{ineq4:thm1:CharWQC:Analo:prev}
         \widetilde{z}_{i}-x=\frac{\alpha_i}{\alpha_{i+1}}(\widetilde{z}_{i+1}-x),
    \end{eqnarray}
    This implies that $ \widetilde{z}_{i}\in [x, \widetilde{z}_{i+1}]$ for all $i=0,1,2,\dots, n-1$. By Theorem \ref{thm1:CharWQC}, we get
    \begin{eqnarray}\label{ineq5:thm1:CharWQC:Analo:prev}
        f(\widetilde{z}_{i})\leq f(\widetilde{z}_{i+1})+\frac{\rho}{2}\|\widetilde{z}_{i+1}-\widetilde{z}_{i}\|\|\widetilde{z}_{i}-x\|.
    \end{eqnarray}
    Summing over $i$, we get
    \begin{eqnarray}\label{ineq5:thm1:CharWQC:Analo:prev}
        f(z)\leq f(y)+\frac{\rho}{2}\sum_{i=0}^{n-1}\|y-z\|\frac{\|\widetilde{z}_{i}-x\|}{n}.
    \end{eqnarray}   
    This implies that 
    \begin{eqnarray}\label{ineq6:thm1:CharWQC:Analo:prev}
        f(z)\leq f(y)+\frac{\rho}{2}\|y-z\|\sum_{i=0}^{n-1}\frac{1}{n}(\|z-x\|+\frac{i}{n}\|y-z\|).
    \end{eqnarray}   
    Therefore, by letting $n\rightarrow\infty$, we obtain
    \begin{eqnarray}\label{ineq7:thm1:CharWQC:Analo:prev}
        f(z)\leq f(y)+\frac{\rho}{2}\|y-z\|\int_{0}^1(\|z-x\|+s\|y-z\|)ds.
    \end{eqnarray}   
    This implies that 
    \begin{eqnarray}\label{ineq8:thm1:CharWQC:Analo:prev}
        f(z)\leq f(y)+\frac{\rho}{2}\|y-z\|(\|z-x\|+\frac{1}{2}\|y-z\|).
    \end{eqnarray} 
    Thus
    \begin{eqnarray}\label{ineq9:thm1:CharWQC:Analo:prev}
        f(z)\leq f(y)+\frac{\rho}{4}(1-\lambda^2)\|y-x\|^2.
    \end{eqnarray} 
\end{proof}

The following example illustrates Theorem \ref{thm1:CharWQC:Analo:prev}.

\begin{example}\label{example1:thm1:CharWQC:Analo:prev}
    Let $C=\mathbb{R}$. Define a function $f:\mathbb{R}\rightarrow\mathbb{R}$ by
    \begin{eqnarray}\label{ineq1:example1:thm1:CharWQC:Analo:prev}
      f(x)=
      \begin{cases}
          0,\;&x<\frac{1}{2}\\
          \frac{1}{2},\;&x\geq\frac{1}{2}.
      \end{cases}
    \end{eqnarray}
    Let $\rho>0$ be a constant. Since $f$ is monotonically increasing, it is quasi-convex on $\mathbb{R}$. Therefore, it is $\rho$-weakly quasi-convex for any $\rho>0$.

    Assume that $f(x)\leq f(\omega)$ for all $\omega\in (z,y]$, where $z=\lambda y+(1-\lambda)x$ for $\lambda\in (0,1]$. If $z$ lies on the left of $y$ then $f(z)\leq f(y)+(\rho/4)(1-\lambda^2)\|y-x\|^2.$ When $y<1/2\leq z$ does not imply that $f(x)\nleq f(\omega)$ for all $\omega\in (z,y]$. 
\end{example}

However, the converse of the above theorem does not hold for a discontinuous function. That is, let $C$ be a nonempty convex subset of $\mathbb{R}^n$ and $f:C\rightarrow\mathbb{R}$ be a discontinuous function. Let $x,y\in C$ and $z=\lambda y+(1-\lambda)x$ for $\lambda\in (0,1]$. Assume that $C^{x,y}_z:=\{\omega\in(z,y]\;\vert\; f(x)\leq f(\omega)\}=(z,y]$ and
    \begin{equation}\label{ineq1:thm1:CharWQC:Analo:prev}
        \begin{aligned}
          f(z) \leq f(y)+\frac{\rho}{4}(1-\lambda^2)\|y-x\|^2,
        \end{aligned}
    \end{equation}
for some $\rho>0.$ Then $f$ may not be a $\rho$-weakly quasi-convex on $C$.

\begin{example}\label{example1:thm1:CharWQC:Analo:prev}
    Let $C=(-\frac{\pi}{2},\frac{\pi}{2})$. Define a function $f:C\rightarrow\mathbb{R}$ by
    \begin{eqnarray}\label{ineq1:example1:thm1:CharWQC:Analo:prev}
      f(x)=
      \begin{cases}
          1,\;&x=0\\
          \sin \vert x\vert,\;&x\neq 0.
      \end{cases}
    \end{eqnarray}
    Let $\rho>0$ be a constant.
    Let $\lambda=1/2$. Take $\varepsilon>0$ such that $x=-\varepsilon$, $y=\varepsilon$ and $x,y\in C$. Then $f(x)=f(y)=\sin\varepsilon$.
    It is obvious that $\lim_{\varepsilon\rightarrow 0^+}(\sin\varepsilon+\frac{\rho\varepsilon^2}{2})=0$, and hence for sufficiently small $\varepsilon>0$, $f(0)=1>\sin\varepsilon+\frac{\rho\varepsilon^2}{2}$. Therefore, $f$ is not $\rho$-weakly quasi-convex on $C.$ 
    
    Assume that $f(x)\leq f(\omega)$ for all $\omega\in (z,y]$, where $z=\lambda y+(1-\lambda)x$ for $\lambda\in (0,1]$. If $x=0$, then $f(x)=f(\omega)=1$ for all $\omega\in (z,y]$ implies $f(y)=1$, and therefore, $y=0$. The inequality  $f(z)\leq f(y)+(\rho/4)(1-\lambda^2)\|y-x\|^2$ holds trivially. Now we claim that $(z,y]$ does not contain $0$. Suppose to the contrary that it contains $0$ then there exists an $\omega$ sufficiently close to $0$, $f(\omega)<f(x)$ since $f(x)>0$ and $\lim_{\omega\rightarrow 0}f(\omega)=\lim_{\omega\rightarrow 0} \sin\vert\omega\vert=0$. Additionally, this implies that $z$ and $y$ lie on the same side of $0$. Now we claim that $f(z)\leq f(y)$. In contrast, assume that $f(z)> f(y)$, which implies that $\sin \vert z\vert>\sin \vert y\vert$. Thus $\vert z\vert>\vert y\vert$. 
    Consider two cases. If $z,y> 0$, then $z>y>0$. This implies that $x>z>y>0$. Similarly, when $z,y<0$, then $x<z<y<0$. Therefore, $\vert z\vert>\vert y\vert$ and $z\in (x,y]$ implies $\vert x\vert>\vert z\vert$. Hence $f(x)>f(z)$. We can now choose an $\omega\in (z,y]$ close to $z$ such that $f(\omega)<f(x)$. Another case when $z=0$ and $y\neq 0$ is not possible. Indeed, if $z=0$, then $\omega= ty$ for $0<t\leq 1$. Therefore, for sufficiently small $t$, $f(w)<f(x)=\sin\vert x\vert$, a contradiction to our hypothesis. Finally, $f(x)\leq f(\omega)$ for all $\omega\in (z,y]$ implies $f(z)\leq f(y)+(\rho/4)(1-\lambda^2)\|y-x\|^2.$
\end{example}

\begin{remark}\label{rmk1:example1:thm1:CharWQC:Analo:prev}

    

    In Theorem \ref{thm1:CharWQC:Analo:prev}, if $f$ is $\rho$-weakly quasi-convex on $C$ that does not imply that for every $x,y\in C$ and every $z=\lambda y+(1-\lambda)x$, $\lambda\in (0,1]$
    \begin{equation}\label{rmk1:ineq1:thm1:CharWQC:Analo:prev}
        \begin{aligned}
          f(x)\leq f(z) \Rightarrow f(z) \leq f(y)+\frac{\rho}{4}(1-\lambda^2)\|y-x\|^2.
        \end{aligned}
    \end{equation}
    For example consider $C=\mathbb{R}$ and $f(x)=-x^2/2$. Then it is easy to check that $f$ is $1$-weakly quasi-convex on $\mathbb{R}$. If we take $x=-1/4$, $y=3/4$, and $\lambda=1/2$ implies $f(x)=f(z)$. However, $f(z) > f(y)+\frac{\rho}{4}(1-\lambda^2)\|y-x\|^2$. 
    Therefore, the assumption in Theorem \ref{thm1:CharWQC:Analo:prev} that $f(x)\leq f(\omega)$ for every $\omega\in(z,y]$ cannot be replaced by the assumption that $f(x)\leq f(z)$.
\end{remark}
The following result is the converse of the above Theorem \ref{thm1:CharWQC:Analo:prev}. We show that the converse of Theorem \ref{thm1:CharWQC:Analo:prev} with an adjusted weak quasi-convexity modulus holds under the continuity assumption of $f$ on the whole domain $C$.

\begin{theorem}\label{thm2:CharWQC:Analo:prev_conv}
    Let $C$ be a nonempty convex subset of $\mathbb{R}^n$ and $f:C\rightarrow\mathbb{R}$ be a continuous function. Let $x,y\in C$ and $z=\lambda y+(1-\lambda)x$, $\lambda\in (0,1]$. Assume that $C^{x,y}_z:=\{\omega\in(z,y]\;\vert\; f(x)\leq f(\omega)\}=(z,y]$ and 
    \begin{equation}\label{ineq1:thm1:CharWQC:Analo:prev}
        \begin{aligned}
          f(z) \leq f(y)+\frac{\rho}{4}(1-\lambda^2)\|y-x\|^2.
        \end{aligned}
    \end{equation}
    Then $f$ is $3\rho/2$-weakly quasi-convex on $C$.
\end{theorem}

\begin{proof}
    Fix arbitrary $x,y\in C$ and $\lambda\in (0,1)$. Define $M:=\max\{f(x), f(y)\}$ and $\varphi(\lambda):=f(\lambda y+(1-\lambda)x)$.
    If $\varphi(\lambda)\leq M$, then the result is straightforward. Otherwise,  $\varphi(\lambda)> M$. Since $f$ is continuous, $\varphi$ is so. Define
    \begin{eqnarray}
        \mathcal{T}_M:=\{t\in (0,1)\;\vert\;\varphi(t)>M\}.
    \end{eqnarray}
    Since $\varphi$ is continuous, $\mathcal{T}_M$ is open in $(0,1)$. It is also obvious that $\mathcal{T}_M$ is nonempty.
    Using the well known fact that every connected subset of $\mathbb{R}$ is an interval, let $(a,b)\subset (0,1)$ be the connected component of $\mathcal{T}_M$ containing $\lambda$. This implies that $a<\lambda<b$.
    \begin{figure}[htbp]
    \centering
    \includegraphics[width=0.5\textwidth]{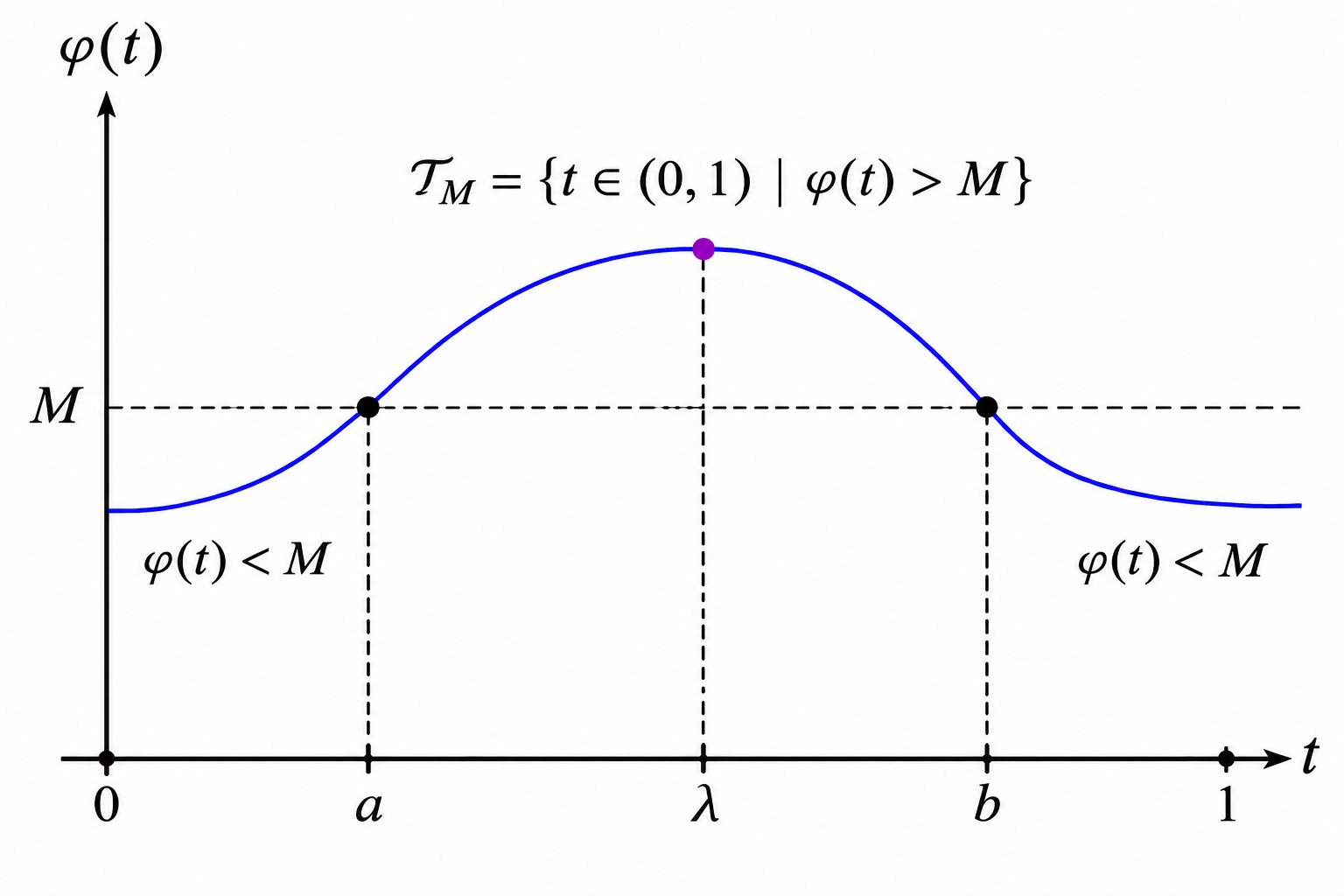}
    \caption{Geometry of $\mathcal{T}_M$}
    \label{fig:geometry of tau_M}
    \end{figure}
    We claim that $\varphi(a)=\varphi(b)=M$. Indeed, let $\{t_n\}$ be a sequence in $\mathcal{T}_M$ such that $t_n\downarrow a$. Then, by the continuity of $\varphi,$ we get, $\varphi(a)\geq M$. Now again claim that $\varphi(a)>M$ cannot hold. Contrary, if true, then by the continuity of $\varphi$ we can get $\varepsilon>0$ sufficiently small such that $\varphi(t)>M$ for all $t\in (a-\varepsilon, a+\varepsilon)$. This contradicts the fact that $(a,b)$ is the maximal connected component of $\mathcal{T}_M$ such that $\varphi(t)>M$. Therefore, $\varphi(a)=M$. Similarly, $\varphi(b)=M$.
    Now consider 
    \begin{eqnarray}
        \tau=\frac{\lambda-a}{b-a}\in (0,1),\; x_a=a y+(1-a)x,\;x_b=b y+(1-b)x.
    \end{eqnarray}
    Then $z_{\lambda}=\tau x_b+(1-\tau)x_a$. Now we claim that $f(x_a)=M\leq f(\omega)$ for all $\omega\in (z_{\lambda},x_b]$.
    Indeed, let $\omega\in(z_{\lambda},x_b]$. Then
    \begin{eqnarray}
    f(x_b)=\varphi(b)=M\leq f(\omega).
    \end{eqnarray}
    When $t=b$ and $\omega=x_b$, then $f(\omega)=f(x_b)=\varphi(b)=M$, and when $\lambda<t<b$, then $f(\omega)>M$.
    Let $\gamma(t)=ty+(1-t)x$ for $t\in (0,1)$. Then for some $\delta\in (0,1)$
    \begin{eqnarray}
    \omega&=&\delta x_b +(1-\delta)z_{\lambda} \nonumber\\
    &=&\delta\gamma(b)+(1-\delta)\gamma(\lambda)\\
    &=&\big((1-\delta)\lambda+\delta b\big)y
    +\big(1-(1-\delta)\lambda-\delta b\big)x\nonumber\\
    &=& \gamma((1-\delta)\lambda+\delta b)\nonumber.
    \end{eqnarray}
    Let $t=(1-\delta)\lambda+\delta b$. Then it is observed that 
    \begin{eqnarray}
        (z_{\lambda},x_b]=\{\gamma(t)\;\vert\; \lambda<t\leq b\}.
    \end{eqnarray}
    Therefore, $f(\omega)\geq M$ for all $\omega\in(z_{\lambda},x_b]$.
    This implies that
    \begin{eqnarray}
    f(x_a)=M\leq f(\omega),
    \end{eqnarray}
    for all $\omega\in(z_{\lambda},x_b]$.
    Therefore, by our hypothesis,
    \begin{eqnarray}
        f(z_{\lambda}) \leq f(x_b)+\frac{\rho}{4}(1-\tau^2)\|x_b-x_a\|^2.
    \end{eqnarray}
    This implies that
    \begin{eqnarray}
        f(z_{\lambda}) \leq M+\frac{\rho}{4}(1-\tau^2)\|x_b-x_a\|^2.
    \end{eqnarray} 
    Now reversing two points, i.e., $z_{\lambda}=(1-\tau)x_a+(1-(1-\tau)) x_b$, we get
    \begin{eqnarray}
        f(z_{\lambda}) \leq M+\frac{\rho}{4}(1-(1-\tau)^2)\|x_b-x_a\|^2.
    \end{eqnarray} 
    Therefore,
    \begin{eqnarray}
        f(z_{\lambda}) \leq M+\frac{\rho}{4}\min\{(1-\tau^2), (1-(1-\tau)^2)\}\|x_b-x_a\|^2.
    \end{eqnarray}  
    It is clear that for $0<\tau< 1$, $\min\{(1-\tau^2), (1-(1-\tau)^2)\}\leq 3\tau(1-\tau)$. Therefore,
    \begin{eqnarray}
        f(\lambda y+(1-\lambda)x) \leq\max\{f(y),f(x)\}+\frac{3\rho}{4}\tau(1-\tau)\|x_b-x_a\|^2.
    \end{eqnarray}
    
    Since
    \begin{eqnarray}
    \tau=\frac{\lambda-a}{b-a},
    \qquad
    1-\tau=\frac{b-\lambda}{b-a},
    \end{eqnarray}
    and
    \begin{eqnarray}
    x_b-x_a=(b-a)(y-x),
    \end{eqnarray}
    we obtain
    \begin{eqnarray}
    \begin{aligned}
    \tau(1-\tau)\|x_b-x_a\|^2
    &=
    \frac{\lambda-a}{b-a}
    \frac{b-\lambda}{b-a}
    (b-a)^2\|y-x\|^2 \\
    &=
    (\lambda-a)(b-\lambda)\|y-x\|^2.
    \end{aligned}
    \end{eqnarray}
    Since $0\leq a<\lambda<b\leq 1$, we have $\lambda-a\le \lambda$ and $b-\lambda\le 1-\lambda$.
    Therefore, $(\lambda-a)(b-\lambda)\leq\lambda(1-\lambda).$
    This implies that
    \begin{eqnarray}
    \tau(1-\tau)\|x_b-x_a\|^2
    \le
    \lambda(1-\lambda)\|y-x\|^2.
    \end{eqnarray}
    Thus
    \begin{eqnarray}
        f(\lambda y+(1-\lambda)x) \leq \max\{f(y),f(x)\}+\frac{3\rho}{4}\lambda(1-\lambda)\|y-x\|^2.
    \end{eqnarray}
    Since $\lambda\in (0,1)$ and $x,y\in C$ are arbitrary, the proof is complete.
\end{proof}

\subsection{First-order conditions and generalized monotonicity}

\begin{proposition}\label{pro1:FirstOrderChar}
 Let $C$ be a nonempty, open and convex subset of $\mathbb{R}^n$. Suppose $f:C\rightarrow\mathbb{R}$ is differentiable and $\rho$-weakly quasi-convex on $C$. Then for any $x,y\in C$, we have
\begin{equation}\label{ineq0:pro1:FirstOrderChar}
        \begin{aligned}
            f(y)\leq f(x) ~ \Longrightarrow ~ \left\langle \nabla f(x), y-x\right\rangle \leq \frac{\rho}{2}\|y-x\|^2.
        \end{aligned}
    \end{equation}
\end{proposition}

\begin{proof}
 Let $\lambda\in (0,1)$ and $x,y\in C$ be such that $f(y)\leq f(x)$. Since $C$ is convex, $\lambda y+(1-\lambda)x\in C$. By the definition of weakly quasi-convexity of $f$, we get   
\begin{equation}\label{ineq2:pro1:FirstOrderChar}
     \begin{aligned}
         f(\lambda y+(1-\lambda)x)\leq f(x)+\frac{\rho}{2}\lambda(1-\lambda)\|y-x\|^2.
     \end{aligned}
 \end{equation}
Therefore,
\begin{equation}\label{ineq2:2:pro1:FirstOrderChar}
     \begin{aligned}
           \frac{f(x+\lambda (y-x))- f(x)}{\lambda}\leq \frac{\rho}{2}(1-\lambda)\|y-x\|^2.
     \end{aligned}
 \end{equation}
Taking the limit as $\lambda\downarrow 0$ in the above inequality, we get
  \begin{equation}\label{ineq3:pro1:FirstOrderChar}
     \begin{aligned}
       \left\langle \nabla f(x), y-x\right\rangle\leq \frac{\rho}{2}\|y-x\|^2,
     \end{aligned}
 \end{equation}
\end{proof}

The following example illustrates Proposition \ref{pro1:FirstOrderChar}.
\begin{example}
    Define a function $f:\mathbb{R}\rightarrow\mathbb{R}$ by
   \begin{equation}
        f(x)=
    \begin{cases}
        -x, & x\leq 1\\
        -x-\frac{(x-1)^2}{2}, &x>1.
    \end{cases}
    \end{equation}
    It is a quasi-convex function on $\mathbb{R}$ since it is monotonically decreasing. Therefore, $f$ is $\rho$-weakly quasi-convex for any $\rho>0$. Also, it is differentiable on $\mathbb{R}$ and its derivative is 
     \begin{equation}
         f^\prime(x)=
    \begin{cases}
        -1, & x\leq 1\\
        -x, &x>1.
    \end{cases}
    \end{equation}
     Any $x\geq y>1$ implies that $f(x)\leq f(y)$. Therefore, $f^\prime (y)(x-y)-(\rho/2)(y-x)^2=-y(x-y)-(\rho/2)(y-x)^2\leq 0$. If $x\geq y$ with $x>1$ and $y\leq 1$ implies $f(x)\leq f(y)$, then $f^\prime (y)(x-y)-(\rho/2)(y-x)^2=-(x-y)-(\rho/2)(y-x)^2\leq 0$. On the other hand, if $1\geq x\geq y$ implies $f(x)\leq f(y)$, then $f^\prime (y)(x-y)-(\rho/2)(y-x)^2=-(x-y)-(\rho/2)(y-x)^2\leq 0$.
\end{example}

The following result is the converse of Proposition \ref{pro1:FirstOrderChar}.

\begin{proposition}\label{pro2:FirstOrderChar}
    Let $C$ be a nonempty, open and convex subset of $\mathbb{R}^n$. Suppose $f:C\rightarrow\mathbb{R}$ is differentiable on $C$. Suppose, there exists a $\rho>0$ such that, for any $x,y\in C$,
\begin{equation}\label{ineq0:pro2:FirstOrderChar}
        \begin{aligned}
            f(y)\leq f(x) ~ \Longrightarrow ~ \left\langle\nabla f(x), y-x\right\rangle\leq\frac{\rho}{2}\|y-x\|^2.
        \end{aligned}
    \end{equation}
    Then $f$ is $(3\rho/2)$-weakly quasi-convex on $C$.
\end{proposition}

\begin{proof}
Let $x,y\in C$ and $z_\lambda=\lambda y+(1-\lambda)x$ for $\lambda\in [0.1]$. Define $g(\lambda)=f(z_\lambda)=f(\lambda y+(1-\lambda)x)$ and 
\begin{eqnarray}\label{ineq0.1:pro2:FirstOrderChar}
    D=\{\lambda\in [0,1]\;\vert\; g(\lambda)>M\},
\end{eqnarray}
where $M=\max\{g(0), g(1)\}=\max\{f(y), f(x)\}$.

Now, if $g(\lambda)\leq M$, then the result is obvious. That is, the result is immediate when $D$ is empty. Otherwise, assume that there exists a $\lambda\in [0,1]$ such that $g(\lambda)>M$. Therefore, $D$ is nonempty. Since $f$ is differentiable on $C$, $g$ is continuous on $[0,1]$. Therefore, $D$ is open in $[0,1]$. Furthermore, $g$ is differentiable on $(0,1)$ and its derivative is $g^\prime(\lambda)=\left\langle\nabla f(z_\lambda), y-x\right\rangle$. Let $(\alpha,\beta)\subset (0,1)$ be a maximal connected component of $D$ containing $\lambda$. Therefore, $g(\mu)>M$ for all $\mu\in (\alpha,\beta)$. This implies that $g(\alpha)\leq M$ and $g(\beta)\leq M$. Furthermore, this implies that $g(\alpha)\leq M<g(\mu)$ and $g(\beta)\leq M< g(\mu)$. 

Now if $g(\alpha)\leq M<g(\mu)$, then $f(z_\alpha)<f(z_\mu)$. Therefore, by our assumption, we get,
\begin{eqnarray}\label{ineq1:pro2:FirstOrderChar}
    \left\langle\nabla f(z_\mu), z_\alpha-z_\mu\right\rangle\leq\frac{\rho}{2}\|z_\alpha-z_\mu\|^2.
\end{eqnarray}
Using the fact $\alpha-\mu<0$ from the above inequality \eqref{ineq1:pro2:FirstOrderChar}, we get
\begin{eqnarray}\label{ineq2:pro2:FirstOrderChar}
    \left\langle\nabla f(z_\mu), y-x\right\rangle\geq-\frac{\rho}{2}(\mu-\alpha)\|y-x\|^2.
\end{eqnarray}
Therefore,
\begin{eqnarray}\label{ineq3:pro2:FirstOrderChar}
    g^\prime(\mu)\geq-\frac{\rho}{2}(\mu-\alpha)\|y-x\|^2.
\end{eqnarray}

Similarly, if $g(\beta)\leq M<g(\mu)$, then $f(z_\beta)<f(z_\mu)$. Therefore, by our assumption, we get
\begin{eqnarray}\label{ineq4:pro2:FirstOrderChar}
    \left\langle\nabla f(z_\mu), z_\beta-z_\mu\right\rangle\leq\frac{\rho}{2}\|z_\beta-z_\mu\|^2.
\end{eqnarray}
Thus, using the fact $\mu-\beta<0$, we obtain
\begin{eqnarray}\label{ineq5:pro2:FirstOrderChar}
    g^\prime(\mu)\leq\frac{\rho}{2}(\beta-\mu)\|y-x\|^2.
\end{eqnarray}

From \eqref{ineq3:pro2:FirstOrderChar} and \eqref{ineq5:pro2:FirstOrderChar}, we get
\begin{eqnarray}\label{ineq6:pro2:FirstOrderChar}
   -\frac{\rho}{2}(\mu-\alpha)\|y-x\|^2\leq g^\prime(\mu)\leq\frac{\rho}{2}(\beta-\mu)\|y-x\|^2.
\end{eqnarray}

Fix a $\delta$ such that $\alpha<\delta<\beta$. Integrating both sides in \eqref{ineq3:pro2:FirstOrderChar} from $\delta$ to $\beta$, we obtain
\begin{eqnarray}\label{ineq7:pro2:FirstOrderChar}
    g(\delta)\leq M+\frac{\rho}{2}(\beta-\delta)(\delta-\alpha)\|y-x\|^2+\frac{\rho}{4}(\beta-\delta)^2\|y-x\|^2.
\end{eqnarray}

Similarly, integrating both sides in \eqref{ineq5:pro2:FirstOrderChar} from $\alpha$ to $\delta$, we obtain
\begin{eqnarray}\label{ineq8:pro2:FirstOrderChar}
    g(\delta)\leq M+\frac{\rho}{2}(\beta-\delta)(\delta-\alpha)\|y-x\|^2+\frac{\rho}{4}(\delta-\alpha)^2\|y-x\|^2.
\end{eqnarray}

Thus, from \eqref{ineq7:pro2:FirstOrderChar} and \eqref{ineq8:pro2:FirstOrderChar}, we get
\begin{eqnarray}\label{ineq9:pro2:FirstOrderChar}
    g(\delta)\leq M+\frac{\rho}{2}(\beta-\delta)(\delta-\alpha)\|y-x\|^2+\frac{\rho}{4}\min\{(\delta-\alpha)^2, (\beta-\delta)^2\}\|y-x\|^2.
\end{eqnarray}

Using the fact that $\min\{a^2,b^2\}\leq ab$ for any $a,b\geq 0$, from \eqref{ineq9:pro2:FirstOrderChar} we obtain,
\begin{eqnarray}\label{ineq10:pro2:FirstOrderChar}
    g(\delta)\leq M+\frac{3\rho}{4}(\beta-\delta)(\delta-\alpha)\|y-x\|^2.
\end{eqnarray}

Thus, 
\begin{eqnarray}\label{ineq11:pro2:FirstOrderChar}
    g(\delta)\leq M+\frac{3\rho}{4}\delta(1-\delta)\|y-x\|^2.
\end{eqnarray}

Since $\delta$, $x$, and $y$ are arbitrary,
\begin{eqnarray}\label{ineq12:pro2:FirstOrderChar}
    f(\delta y+(1-\delta)x)\leq\max\{f(y), f(x)\}+\frac{3\rho}{4}\delta(1-\delta)\|y-x\|^2,
\end{eqnarray}
for all $\delta\in (\alpha,\beta)\subset (0,1)$ and $x,y\in C$.

Thus, for any $x,y\in C$ and $\lambda\in [0,1]$
\begin{eqnarray}\label{ineq12:pro2:FirstOrderChar}
    f(\lambda y+(1-\lambda)x)\leq\max\{f(y), f(x)\}+\frac{3\rho}{4}\lambda(1-\lambda)\|y-x\|^2.
\end{eqnarray}
\end{proof}

\begin{remark}
   The constructive nature of the proof of the above theorem naturally raises the question of whether the weak quasi-convexity modulus of $f$ can be improved from $3\rho/2$ to $\rho$.
\end{remark}

The following example illustrates Proposition \ref{pro2:FirstOrderChar}.

\begin{example}\label{example1:pro2:FirstOrderChar}
    Let $C=(0,1)$ and consider the function $f:C\rightarrow\mathbb{R}$ by
    \begin{eqnarray}\label{fun1:example1:pro2:FirstOrderChar}
    f(x)=
    \begin{cases}
        7x(1-x),\;&0<x\leq \frac{1}{2}\\
        \frac{7}{4}-\frac{1}{4}\big(x-\frac{1}{2}\big)^2,\;&\frac{1}{2}<x<1.
    \end{cases}
    \end{eqnarray}

    Note that $f$ is differentiable on $C$ and its derivative is obtained by
    \begin{eqnarray}\label{derifun1:example1:pro2:FirstOrderChar}
    f^\prime(x)=
    \begin{cases}
        14\big(\frac{1}{2}-x\big),\;&0<x\leq \frac{1}{2}\\
        -\frac{1}{2}\big(x-\frac{1}{2}\big),\;&\frac{1}{2}<x<1.
    \end{cases}
    \end{eqnarray}

    \textbf{Case-I}: Assume that $x,y\in C$ and $x,y\leq 1/2$. Then $f(y)\leq f(x)$ implies $y\leq x$. Therefore,
    \begin{eqnarray}\label{ineq1:example1:pro2:FirstOrderChar}
         f'(x)(y-x)-\frac{\rho}{2}(y-x)^2=(y-x)\left(14\left(\frac{1}{2}-x\right)-\frac{\rho}{2}(y-x)\right)\leq 0
    \end{eqnarray}
    for any $\rho>0$.

    \textbf{Case-II}: Suppose that $x,y\in C$ and $x,y> 1/2$. Then $f(y)\leq f(x)$ implies $y\geq x$. Therefore,
    \begin{eqnarray}\label{ineq2:example1:pro2:FirstOrderChar}
         f'(x)(y-x)-\frac{\rho}{2}(y-x)^2=(y-x)\left(\left(\frac{1}{2}-x\right)-\frac{\rho}{2}(y-x)\right)\leq 0
    \end{eqnarray}
    for any $\rho>0$.

    \textbf{Case-III}: Assume that $x,y\in C$ and $y\leq 1/2<x$. Let $y=1/2-\varepsilon_1$ and $x=1/2+\varepsilon_2$, where $\varepsilon_1\geq 0$ and $\varepsilon_2>0$.  Then $f(y)\leq f(x)$ implies $28\varepsilon_1(1-2\varepsilon_1)\leq 7-\varepsilon_2^2$, which further implies $\varepsilon_1\geq \varepsilon_2/\sqrt{28}$. Therefore, 
    \begin{eqnarray}\label{ineq3:example1:pro2:FirstOrderChar}
     \frac{\varepsilon_2}{\varepsilon_1+\varepsilon_2}\leq \frac{\varepsilon_2}{\frac{\varepsilon_2}{\sqrt{28}}+\varepsilon_2}.
    \end{eqnarray}
    
    Now assume that 
    \begin{eqnarray}\label{ineq4:example1:pro2:FirstOrderChar}
         \frac{\rho}{2}(y-x)^2-f'(x)(y-x)\geq 0 ~\Longleftrightarrow~  \frac{\rho}{2}(\varepsilon_1+\varepsilon_2)^2-\frac{1}{2}\varepsilon_2(\varepsilon_1+\varepsilon_2)\geq 0.
    \end{eqnarray}
    This implies that 
    \begin{eqnarray}\label{ineq5:example1:pro2:FirstOrderChar}
        \rho\geq \frac{\varepsilon_2}{\varepsilon_1+\varepsilon_2}. 
    \end{eqnarray}
    From \eqref{ineq3:example1:pro2:FirstOrderChar} and \eqref{ineq4:example1:pro2:FirstOrderChar}, we conclude that \eqref{ineq4:example1:pro2:FirstOrderChar} holds for any $\rho\geq \sqrt{28}/(1+\sqrt{28})$.

    \textbf{Case-IV}: Let $x,y\in C$ and $x\leq 1/2<y$. Let $x=1/2-\varepsilon_1$ and $y=1/2+\varepsilon_2$, where $\varepsilon_1\geq 0$ and $\varepsilon_2>0$.  Then $f(y)\leq f(x)$ implies 
    \begin{eqnarray}\label{ineq6:example1:pro2:FirstOrderChar}
        \frac{7}{4}-\frac{\varepsilon_2^2}{4}\leq 7\left(\frac{1}{4}-\varepsilon_1^2\right),
    \end{eqnarray}
    which further implies $\varepsilon_2\geq \sqrt{28}\varepsilon_1.$
    
    Now assume that 
    \begin{eqnarray}\label{ineq7:example1:pro2:FirstOrderChar}
         \frac{\rho}{2}(y-x)^2-f'(x)(y-x)\geq 0 ~\Longleftrightarrow~  \frac{\rho}{2}(\varepsilon_1+\varepsilon_2)^2-14\varepsilon_1(\varepsilon_1+\varepsilon_2)\geq 0.
    \end{eqnarray}
    This implies that 
    \begin{eqnarray}\label{ineq8:example1:pro2:FirstOrderChar}
        \rho\geq \frac{28\varepsilon_1}{\varepsilon_1+\varepsilon_2}. 
    \end{eqnarray}
    Since 
    \begin{eqnarray}\label{ineq8.1:example1:pro2:FirstOrderChar}
        \frac{28\varepsilon_1}{\varepsilon_1+\varepsilon_2}\leq\frac{28}{1+\sqrt{28}},
    \end{eqnarray}
    the inequality in \eqref{ineq7:example1:pro2:FirstOrderChar} holds for any $\rho\geq 28/(1+\sqrt{28})$.

   Combining all the cases we conclude that for $\rho\geq 28/(1+\sqrt{28})$ and for any $x,y\in C$ with
   \begin{eqnarray}\label{ineq9:example1:pro2:FirstOrderChar}
       f(y)\leq f(x)\implies f'(x)(y-x)\leq \frac{\rho}{2}(y-x)^2.
   \end{eqnarray}
   Therefore, \eqref{ineq0:pro2:FirstOrderChar} is satisfied. Hence, by Proposition \ref{pro2:FirstOrderChar}, the function $f$ is $3\rho/2$-weakly quasi-convex on $C$. 
\end{example}

\subsubsection{Weak quasi-monotonicity}

\begin{definition}\label{Def1:WeakQuasiMono}
    Let $C$ be a subset of $\mathbb{R}^n$. A mapping $g:C\rightarrow\mathbb{R}^n$ is said to be $\rho$-weakly quasi-monotone on $C$ for some $\rho>0$ if for any $x,y\in C$,
    \begin{equation}\label{gen:mon}
     \left\langle g(x), y-x \right\rangle - \frac{\rho}{2} \|y-x\|^2 > 0 ~ \Longrightarrow ~ \left\langle g(y), y-x \right\rangle + \frac{\rho}{2} \|y-x\|^2 \geq 0.
    \end{equation}
\end{definition}

\begin{remark}\label{rmk1:Lem1:FirstOrderChar}
\begin{enumerate}
    \item[(1)] Every hypomonotone mapping (see, \cite[Definition 2, on $\mathbb{R}^n$]{KALO2025}), i.e., for any $x,y\in\mathbb{R}^n$ and some $\eta>0$ such that
    \begin{eqnarray}
        \left\langle G(y)-G(x),y-x\right\rangle\geq-\eta\|y-x\|^2,
    \end{eqnarray}
    where $G:\mathbb{R}^n\to\mathbb{R}^n$, is $\rho$-weakly quasi-monotone on $\mathbb{R}^n$ with the same modulus. However, the converse is not true. For example, take $f(x)=-\exp(1+x^2)$ for $x\in \mathbb{R}$. Since $f$ is quasi-monotone, it is weakly quasi-monotone with any modulus $\rho>0.$ However, it is not hypomonotone on $\mathbb{R}$ with modulus $\rho$. Indeed, contrary assume that for any $x,y\in\mathbb{R}$ and some $\rho>0$ such that
    \begin{eqnarray}
        \left\langle f(y)-f(x),y-x\right\rangle\geq-\rho\|y-x\|^2.
    \end{eqnarray}
    Since $f$ is differentiable, taking the limit $y\to x$ in the above inequality we get $f^\prime (x)\geq-\rho$ for all $x\in\mathbb{R}$. This implies that $2x\exp(1+x^2)\leq \rho$ for all $x\in\mathbb{R}$, a contradiction.
    
    \item[(2)]  If we take $\rho=0$, then Definition \ref{Def1:WeakQuasiMono} reduces to the standard definition for quasi-monotonicity of $f$ (see, for instance, \cite[Definition 5.1]{SKAR1990}). In general, every quasi-monotone mapping is weakly quasi-monotone.
\end{enumerate}
\end{remark}

In the next example, we show that a weakly quasi-monotone mapping may not be quasi-monotone. Furthermore, this gives us the guarantee that the class of weakly quasi-monotone mappings is substantially broader than the class of quasi-monotone mappings. 

\begin{example}\label{example1:WeakQuasiMonoNotQM}
 Define a function $f:\mathbb{R} \rightarrow \mathbb{R}$ by
\[
f(x)=
\left\{
\begin{array}{ll}
 x^{2},  & x \neq  0\\
 -1,    &  x = 0.
\end{array}
\right.
\]

    We claim that $f$ is not quasi-monotone but it is weakly quasi-monotone with modulus $2$.
    Indeed, take $x=0$ and $y=-1$. Then $f(y)(y-x)=-1<0$, however $f(x)(y-x)=1>0$. Therefore, $f$ is not quasi-monotone.
    To establish its weak quasi-monotonicity with modulus $\rho=2$, consider the following cases:

    \textbf{Case-I:} Suppose $x\neq 0$. Then $f(x)(y-x)-\frac{\rho}{2}(y-x)^2=(y-x)(x^2+x-y)>0$ implies that $y>x$ and $x^2+x-y>0$ or $y<x$ and $x^2+x-y<0$. If $y<x$ and $x^2+x-y<0$, then $x^2<0$, a contradiction. 

    Now consider the valid case $y>x$ and $x^2+x-y>0$. In this case, we will have the following scenarios:
    \begin{enumerate}
        \item[(i)] if $y\neq 0$, then
    \begin{eqnarray}
        \begin{aligned}
            f(y)(y-x)+\frac{\rho}{2}(y-x)^2=(y-x)(y^2+y-x)\geq 0.
        \end{aligned}
    \end{eqnarray}

        \item[(ii)] if $y=0$, then 
    \begin{eqnarray}
        \begin{aligned}
            f(y)(y-x)+\frac{\rho}{2}(y-x)^2=x+x^2\geq 0.
        \end{aligned}
    \end{eqnarray}
    \end{enumerate}

    \textbf{Case-II:} Suppose $x=0$ and $f(x)(y-x)-\frac{\rho}{2}(y-x)^2=-y(y+1)>0$. We will have the following scenarios: if $y\neq 0$, 
    \begin{eqnarray}
        \begin{aligned}
            f(y)(y-x)+\frac{\rho}{2}(y-x)^2=y^2(y+1).
        \end{aligned}
    \end{eqnarray}
     If $-1<y<0$, then $y^2(y+1)\geq 0.$ Otherwise, contradict the condition $y^2+y<0$.
\end{example}

\begin{proposition}\label{pro1:WeakQuasiMono}
 Let $C$ be a nonempty subset of $\mathbb{R}^n$, and let $f: C \rightarrow\mathbb{R}^n$ be a weakly quasi-monotone mapping on $C$ with modulus $\rho > 0$, and $\gamma>0$ be a constant. Then $\gamma f$ is weakly quasi-monotone on $C$ with modulus $\gamma \rho > 0$.
\end{proposition}

\begin{proof}
  Let $x, y\in C$. Assume that $\left\langle \gamma f(x), y-x \right\rangle - \frac{\gamma\rho}{2} \|y-x\|^2>0$.  This implies that $\left\langle f(x), y-x \right\rangle - \frac{\rho}{2} \|y-x\|^2 > 0$. 
    Using the weak quasi-monotonicity of $f$ on $C$, we get
    \begin{eqnarray}\label{ineq2:pro1:WeakQuasiMono}
        \begin{aligned}
 \left\langle  f(y),y-x \right\rangle + \frac{\rho}{2} \|y-x\|^2 \geq 0 ~ \Longleftrightarrow ~ \left\langle  \gamma f(y),y-x\right\rangle + \frac{\gamma \rho}{2} \|y-x\|^2\geq 0, \notag
        \end{aligned}
    \end{eqnarray}
 and the result follows.   
\end{proof}



\begin{lemma}\label{Lem1:FirstOrderChar}
 Let $C$ be an open convex subset of $\mathbb{R}^n$, and let $f: C \rightarrow \mathbb{R}$ be a differentiable and $\rho$-weakly quasi-convex function with $\rho > 0$. Then $\nabla f$ is $\rho$-weakly quasi-monotone on $C$.
\end{lemma}

\begin{proof}
    Let for any $x,y\in C$, the function $f$ satisfies
    \begin{eqnarray}\label{ineq1:Lem1:FirstOrderChar}
        \begin{aligned}
            \left\langle \nabla f(x), y-x\right\rangle-\frac{\rho}{2}\|y-x\|^2> 0.
        \end{aligned}
    \end{eqnarray}

    From the definition of weakly quasi-convexity of $f$ and \eqref{ineq1:Lem1:FirstOrderChar}, we get $f(x)< f(y).$ This implies that
\begin{align*}\label{ineq2:Lem1:FirstOrderChar}
 \left\langle \nabla f(y), x-y\right\rangle-\frac{\rho}{2}\|y-x\|^2\leq 0 ~ \Longleftrightarrow ~ \left\langle \nabla f(y), y-x\right\rangle +\frac{\rho}{2}\|y-x\|^2 \geq 0,
\end{align*}  
and the result follows.
\end{proof}

\begin{remark}
    Observe that our monotonicity introduced in Definition \ref{Def1:WeakQuasiMono} is weaker than the monotonicity produced by \cite[Proposition 8, implication in (19)]{FLAR2025}, which is the following. Let $C$ be a nonempty subset of $\mathbb{R}^n$ and $g:C\rightarrow\mathbb{R}^n$ is a mapping. Then $g$ is $\delta$-strongly quasi-monotone on $C$ if there exists a $\delta>0$ such that
    \begin{eqnarray}\label{ineq1:gen:mon:com:LARA(19)}
         \left\langle g(x), y-x \right\rangle> - \frac{\delta}{2}\|y-x\|^2 \implies \left\langle g(y), y-x \right\rangle \geq \frac{\delta}{2}\|y-x\|^2.
    \end{eqnarray}
    
    Indeed, suppose that $g$ satisfies the implication in \eqref{ineq1:gen:mon:com:LARA(19)}. To show the implication in \eqref{gen:mon}, assume that there exists a $\rho>0$ such that
    \begin{eqnarray}\label{ineq1:gen:mon:com}
     \left\langle g(x), y-x \right\rangle - \frac{\rho}{2} \|y-x\|^2 > 0.
    \end{eqnarray}
    This implies that 
    \begin{eqnarray}\label{ineq2:gen:mon:com}
     \left\langle g(x), y-x \right\rangle>\frac{\rho}{2} \|y-x\|^2 \geq-\frac{\rho}{2}\|y-x\|^2,
    \end{eqnarray}
    which further implies
    \begin{eqnarray}\label{ineq3:gen:mon:com}
     \left\langle g(y), y-x \right\rangle+\frac{\rho}{2}\|y-x\|^2\geq\rho \|y-x\|^2\geq 0.
    \end{eqnarray}
    Hence, the required inequality follows. Now we provide an example to show the negation of its converse.
    Consider $g(x)=-x$ for $x\in\mathbb{R}$ and $\rho=1$. Then 
    \begin{eqnarray}\label{ineq3.1:gen:mon:com}
     -x(y-x)-\frac{1}{2} (y-x)^2>0\implies -y(y-x)+\frac{1}{2} (y-x)^2> 0.
    \end{eqnarray}    
    However, particularly, if we consider $x=-1$ and $y=1$ then $-x(y-x)>-\frac{\delta}{2} (y-x)^2$ implies $1>-\delta$ always holds but $ -y(y-x)\geq\frac{\delta}{2} (y-x)^2$ implies $-1\geq \delta$, which is impossible.
    
    Therefore, in general, the notion of $\rho$-weakly quasi-monotonicity on $C$ is more flexible and general than the notion defined in \eqref{ineq1:gen:mon:com:LARA(19)}, and therefore applies to a broader class of mappings.

\end{remark}


The example below shows that the converse of the above Lemma \ref{Lem1:FirstOrderChar} is not true. 

\begin{example}\label{example1:Lem1:FirstOrderChar}
    Let $C=(0,1)$ and consider the function $f:C\rightarrow\mathbb{R}$ by
    \begin{eqnarray}
        f(x)=
    \begin{cases}
        7x(1-x),\;&0<x\leq \frac{1}{2}\\
        \frac{7}{4}-\frac{1}{4}(x-\frac{1}{2})^2,\;&\frac{1}{2}<x<1.
    \end{cases}
    \end{eqnarray}

    Note that $f$ is differentiable on $C$ and its derivative is obtained by
    \begin{eqnarray}
        f^\prime(x)=
    \begin{cases}
        14(\frac{1}{2}-x),\;&0<x\leq \frac{1}{2}\\
        -\frac{1}{2}(x-\frac{1}{2}),\;&\frac{1}{2}<x<1.
    \end{cases}
    \end{eqnarray}

    Define $g(x)=f^\prime(x)$ for $x\in C$. We claim that $g$ is $1$-weakly quasi-monotone on $C$. Indeed, let
    $g(x)(y-x)-(1/2)(y-x)^2>0$ for $x,y\in C$. This implies that either $y>x$ and $g(x)>(y-x)/2$ or $y<x$ and $g(x)<(y-x)/2$. 
    
    Let $y>x$ and $g(x)>(y-x)/2$ for $x,y\in C$. 
    If $x>1/2$ and $y>x$, then $y>x> 1/2$ but $g(x)<(y-x)/2$. This is impossible.
    If $x<1/2$, $y>1/2$ and $g(x)>(y-x)/2$, then $g(y)=(1/2)(1/2-y)>-(y-x)/2$. Therefore, $(y-x)(g(y)+(y-x)/2)\geq 0$.
    If $x\leq 1/2$, $y\leq 1/2$ and $y>x$, then $g(x)=14(1/2-x)>14(y-x)>(y-x)/2$ and $g(y)=14(1/2-y)\geq 0>-(y-x)/2$. Therefore, $(y-x)(g(y)+(y-x)/2)\geq 0$. 

    \medskip{}
    Similarly, if $y<x$ and $g(x)<(y-x)/2$. 
    If $y>1/2$ and $x>y$, then $x>y> 1/2$. Therefore, $g(x)=-(1/2)(x-1/2)<(y-x)/2<(y-x)/2$ and $g(y)=-(1/2)(y-1/2)<-(y-x)/2$. Finally, $(y-x)(g(y)+(y-x)/2)\geq 0$.
    If $y\leq 1/2$ and $x>1/2$, then $g(x)\geq (y-x)/2$. This is impossible.  
    If $x\leq 1/2$, $y\leq 1/2$ and $x>y$, then $g(x)=14(1/2-x)\geq 0$. However, $g(x)<(y-x)/2<0$, impossible. 
    
    Finally, $g$ is $1$-weakly quasi-monotone on $C$. However, it fails to satisfy 
    \begin{eqnarray*}
        \begin{aligned}
            f(y)\leq f(x)\Rightarrow g(x)(y-x)\leq (y-x)^2/2.
        \end{aligned}
    \end{eqnarray*}
    Indeed, let $x=0.42$ and $y=0.99$. This implies that $f(x)=1.7052$ and $f(y)=1.6899$. Therefore, $f(y)<f(x)$. Now, $g(x)(y-x)=0.6384$ and $(y-x)^2/2=0.1624$. Therefore, $g(x)(y-x)>(y-x)^2/2$.
    Thus by Proposition \ref{pro1:FirstOrderChar}, $f$ is not $1$-weakly quasi-convex on $C$. 
\end{example}

\begin{remark}\label{rmk1:example1:Lem1:FirstOrderChar}
\begin{enumerate}
    \item Example \ref{example1:Lem1:FirstOrderChar} demonstrates that a function can be $1$-weakly quasi-monotone without being $1$-weakly quasi-convex on $C$.

    \item Since the converse of Proposition \ref{pro1:FirstOrderChar} holds with weakly quasi-convex modulus $3\rho/2$, it is natural to ask whether an analogous converse of Lemma \ref{Lem1:FirstOrderChar} also holds.
    Indeed, let $C$ be a nonempty, open and convex subset of $\mathbb{R}^n$, and let the function $f:C\rightarrow\mathbb{R}$ be differentiable, and $\nabla f$ be $\rho$-weakly quasi-monotone on $C$. 
    Does it then follow that the function $f$ is $3\rho/2$-weakly quasi-convex function on $C$? 
    The answer to this question is negative. Consider all the data mentioned in Example \ref{example1:Lem1:FirstOrderChar}. If we take $x=0.42$ and $y=0.99$, the inequality $g(x)(y-x)\leq3(y-x)^2/4$ does not hold.
\end{enumerate}
\end{remark}


\subsubsection{Weak pseudo-monotonicity and strict weak pseudo-monotonicity}\label{WSWSM}

\begin{definition}\label{Def2:WSWSM}
    Let $C$ be a nonempty subset of $\mathbb{R}^n$ and $\rho>0$ be a constant. A mapping $f:C\rightarrow\mathbb{R}^n$ is said to be 
    \begin{enumerate}
        \item[(1)] $\rho$-weakly pseudo-monotone on $C$, if for any $x,y\in C$
        \begin{eqnarray}\label{ineq1:Def2:WSWSM}
        \begin{aligned}
            \left\langle f(x), y-x\right\rangle-\frac{\rho}{2}\|y-x\|^2\geq 0\implies \left\langle f(y), y-x\right\rangle+\frac{\rho}{2}\|y-x\|^2\geq 0.
        \end{aligned}
    \end{eqnarray}
    \item[(2)] $\rho$-strictly weakly pseudo-monotone on $C$, if for any $x, y\in C$ with $x\neq y$
        \begin{eqnarray}\label{ineq2:Def2:WSWSM}
        \begin{aligned}
            \left\langle f(x), y-x\right\rangle-\frac{\rho}{2}\|y-x\|^2\geq 0\implies \left\langle f(y), y-x\right\rangle+\frac{\rho}{2}\|y-x\|^2> 0.
        \end{aligned}
    \end{eqnarray}
    \end{enumerate}
\end{definition} 

\begin{remark}\label{rmk1:Def2:WSWSM}
\begin{enumerate}
    \item[(1)] Every $\rho$-strictly weakly pseudo-monotone mapping on $C$ is $\rho$-weakly pseudo-monotone on $C$.
    \item[(2)] Clearly, every pseudo-monotone mapping on $C$ (see, \cite[Definition 3.1]{SKAR1990}) is $\rho$-weakly pseudo-monotone on $C$, but the converse does not hold. Indeed, consider the function in Example \ref{example1:WeakQuasiMonoNotQM}. Clearly, it is $2$-weakly pseudo-monotone, but not pseudo-monotone on $\mathbb{R}$.
\end{enumerate}
\end{remark}


\begin{proposition}\label{Pro2:WSWSM}
    Let $C$ be a nonempty subset of $\mathbb{R}^n$. A mapping  $f:C\rightarrow\mathbb{R}^n$ is $\rho$-weakly pseudo-monotone on $C$ if and only if
    \begin{eqnarray}\label{ineq1:Pro2:WSWSM}
        \begin{aligned}
            \left\langle f(x), y-x\right\rangle-\frac{\rho}{2}\|y-x\|^2> 0\implies \left\langle f(y), y-x\right\rangle+\frac{\rho}{2}\|y-x\|^2> 0,\;\;\text{for all}\;\;x,y\in C.
        \end{aligned}
    \end{eqnarray}
\end{proposition}

\begin{proof}
    By definition, $f$ is $\rho$-weakly pseudo-monotone iff 
    \begin{eqnarray}\label{ineq2:Pro2:WSWSM}
        \begin{aligned}
            \left\langle f(x), y-x\right\rangle-\frac{\rho}{2}\|y-x\|^2\geq 0\implies \left\langle f(y), y-x\right\rangle+\frac{\rho}{2}\|y-x\|^2\geq 0,\;\;\text{for all}\;\; x,y\in C.
        \end{aligned}
    \end{eqnarray}
    This implies that 
        \begin{eqnarray}\label{ineq2:Pro2:WSWSM}
        \begin{aligned}
            \left\langle f(y), y-x\right\rangle+\frac{\rho}{2}\|y-x\|^2< 0\implies \left\langle f(x), y-x\right\rangle-\frac{\rho}{2}\|y-x\|^2<0,\;\;\text{for all}\;\; x,y\in C.
        \end{aligned}
    \end{eqnarray}
    Thus,
    \begin{eqnarray}\label{ineq3:Pro2:WSWSM}
        \begin{aligned}
            \left\langle f(y), x-y\right\rangle-\frac{\rho}{2}\|y-x\|^2> 0\implies \left\langle f(x), x-y\right\rangle+\frac{\rho}{2}\|y-x\|^2>0,\;\;\text{for all}\;\; x,y\in C.
        \end{aligned}
    \end{eqnarray}
\end{proof}

\begin{remark}\label{rmk1:Pro2:WSWSM}
    Proposition \ref{Pro2:WSWSM}, guarantees that one cannot generalize the definition of $\rho$-weakly pseudo-monotonicity of $f$ on $C$ by changing inequalities in \eqref{ineq1:Def2:WSWSM}.
\end{remark}

\begin{remark}\label{rmk2:Def2:WSWSMWQM}
    Every $\rho$-weakly pseudo-monotone mapping on $C$ is $\rho$-weakly quasi-monotone on $C$. But the converse is not true. We provide an example below.
\end{remark}

\begin{example}\label{exam2:Def2:WSWSMWQM}
    Let $\rho>0$ be fixed and $C=\{0,1\}$ be a subset of $\mathbb{R}$. Define $f:C\rightarrow\mathbb{R}$ by $f(0)=\rho$ and $f(1)=-\rho/2$. It is a $\rho$-weakly quasi-monotone for any $\rho>0$. Indeed, let $x,y\in\mathbb{R}$ and $f(x)(y-x)-\rho/2(y-x)^2> 0$. This implies that $y>x$ and $f(x)> \rho (y-x)/2$, or $y<x$ and $f(x)< \rho (y-x)/2$. The latter condition is not possible. Therefore, we have only $y>x$ and $f(x)> \rho (y-x)/2$. This says that $y=1$ and $x=0$. Therefore, $f(y)(y-x)+\rho(y-x)^2/2=f(1)+\rho/2=0.$ However, it is not a $\rho$-weakly pseudo-monotone function by Proposition \ref{Pro2:WSWSM}.
\end{example}
\subsection{Second-order necessary conditions and related properties}

\begin{proposition}\label{pro1:SecondOrderChar_Dini}
 Let $C$ be a nonempty, open and convex subset of $\mathbb{R}^n$ and $f: C \rightarrow\mathbb{R}$ be a $\rho$-weakly quasi-convex function. 
 Assume that 
 \begin{enumerate}
     \item[(A1)] for $x\in C$ and $d\in\mathbb{R}^n$, the directional derivative $f^\prime(x;d)$ exists and $f$ has second-order Dini derivative on $C$.
     \item[(A2)]  for $x\in C$ and $d\in\mathbb{R}^n$ such that $t\mapsto f(x+td)$ is absolutely continuous on every compact subinterval of
         \begin{eqnarray*}
             I_{x,d}:=\{t\in \mathbb{R}\;\vert\; x+td\in C\}.
         \end{eqnarray*}
 \end{enumerate}
 Then for $x\in C$ and $d\in\mathbb{R}^n$ such that
    \begin{equation}\label{ineq2:pro1:SecondOrderChar_Dini}
        \begin{aligned}
           f^\prime(x,d)=0=f^\prime(x,-d)\implies\max\left\{D^{\prime\prime}f(x;d),D^{\prime\prime}f(x;-d)\right\}+\rho\|d\|^2\geq 0.
        \end{aligned}
    \end{equation}
\end{proposition}

\begin{proof}
    By the Dini's second-order derivative, we get
    \begin{eqnarray}\label{ineq3:pro1:SecondOrderChar_Dini}
        D^{\prime\prime}f(x;d)=\lim_{\delta\to 0^+}\frac{f'(x+\delta d;d)-f'(x;d)}{\delta},
    \end{eqnarray}
    where
    \begin{eqnarray}\label{ineq4:pro1:SecondOrderChar_Dini}
        f'(x;d)=\lim_{t\to 0^+}\frac{f(x+td)-f(x)}{t}.
    \end{eqnarray}

    From \eqref{ineq3:pro1:SecondOrderChar_Dini}, we get
    \begin{eqnarray}\label{ineq5:pro1:SecondOrderChar_Dini}
    \lim_{\delta\to 0^+}\frac{f'(x+\delta d;d)-f'(x;d)-\delta D^{\prime\prime}f(x;d)}{\delta}=0.
    \end{eqnarray}
    
    Hence,
    \begin{eqnarray}\label{ineq6:pro1:SecondOrderChar_Dini}
    f'(x+\delta d;d)=f'(x;d)+\delta D^{\prime\prime}f(x;d)+o_+(\delta),
    \end{eqnarray}
    where 
    \begin{eqnarray}
        \frac{o_+(\delta)}{\delta}\to 0,\;\; \text{as }\delta\to 0^+.
    \end{eqnarray}

    Let $h> 0$ such that $[0,h]\subset I_{x,d}$ and define
    \begin{eqnarray}\label{ineq7:pro1:SecondOrderChar_Dini}
    \varphi(t):=f(x+td),\;\;\text{for all}\;\; t\in[0,h].
    \end{eqnarray}
    Using our hypothesis (A2), by the fundamental theorem of calculus for absolutely continuous functions,
    \begin{eqnarray}\label{ineq8:pro1:SecondOrderChar_Dini}
    \varphi(h)-\varphi(0)=\int_0^h \varphi'(t)\,dt.
    \end{eqnarray}

    Since every absolutely continuous function is differentiable almost
    everywhere, at every point $t$ at which $\varphi$ is differentiable,
    \begin{eqnarray}\label{ineq9:pro1:SecondOrderChar_Dini}
    \begin{aligned}
    \varphi'(t)&=\lim_{s\to0}\frac{\varphi(t+s)-\varphi(t)}{s}\\
    &=\lim_{s\to0}\frac{f(x+td+sd)-f(x+td)}{s}.
    \end{aligned}
    \end{eqnarray}
    
    By our hypothesis (A2)
    \begin{eqnarray}\label{ineq10:pro1:SecondOrderChar_Dini}
    \varphi'(t)=f'(x+td;d)\;\; \text{for a.e. } t.
    \end{eqnarray}
    
    Therefore,
    \begin{eqnarray}\label{ineq11:pro1:SecondOrderChar_Dini}
    f(x+hd)-f(x)=\int_0^h f'(x+td;d)\,dt.
    \end{eqnarray}

    Replacing $\delta$ with $s$ in $\eqref{ineq6:pro1:SecondOrderChar_Dini}$, we get
    \begin{eqnarray}\label{ineq12:pro1:SecondOrderChar_Dini}
    f'(x+sd;d)=f'(x;d)+sD^{\prime\prime}f(x;d)+o_+(s).
    \end{eqnarray}
    Using (A2) and $f^\prime(x;d)=0$, we obtain
    \begin{eqnarray}\label{ineq13:pro1:SecondOrderChar_Dini}
    \begin{aligned}
    f(x+hd)-f(x)&=\int_0^h f'(x+sd;d)\,ds\\
    &=\int_0^h\left(sD^{\prime\prime}f(x;d)+o_+(s)\right)\,ds\\
    &=\frac{h^2}{2}D^{\prime\prime}f(x;d)+o_+(h^2).
    \end{aligned}
    \end{eqnarray}

    We claim that
    \begin{eqnarray}
        \int_0^h o_+(s)\,ds=o_+(h^2).
    \end{eqnarray}

    Indeed, for any $\varepsilon>0$
    \begin{eqnarray}
        \frac{1}{h^2}\Bigg\vert\int_0^h o_+(s)\,ds\Bigg\vert\leq\frac{1}{h^2}\Bigg\vert\int_0^h \varepsilon s\,ds\Bigg\vert=\frac{\varepsilon}{2}.
    \end{eqnarray}
    
    Thus,
    \begin{eqnarray}\label{ineq14:pro1:SecondOrderChar_Dini}
    f(x+hd)=f(x)+\frac{h^2}{2}D^{\prime\prime}f(x;d)+o_+(h^2).
    \end{eqnarray}
    
    Similarly, by our hypotheses we obtain
    \begin{eqnarray}\label{ineq15:pro1:SecondOrderChar_Dini}
    f(x-hd)=f(x)+\frac{h^2}{2}D^{\prime\prime}f(x;-d)+o_{-}(h^2).
    \end{eqnarray}

    Since $f$ is weakly quasi-convex on $C$,
    \begin{eqnarray}\label{ineq16:pro1:SecondOrderChar_Dini}
    f(x)\leq\max\left\{f(x+hd),\,f(x-hd)\right\}+\frac{\rho}{2}h^2\|d\|^2.
    \end{eqnarray}
    
    Therefore,
    \begin{eqnarray}\label{ineq17:pro1:SecondOrderChar_Dini}
    \begin{aligned}
    f(x)&\leq\max\Bigg\{f(x)+\frac{h^2}{2}D^{\prime\prime}f(x;d)+o_+(h^2),\\
    &\hspace{3.5cm}f(x)+\frac{h^2}{2}D^{\prime\prime}f(x;-d)+o_{-}(h^2)\Bigg\}+\frac{\rho}{2}h^2\|d\|^2.
    \end{aligned}
    \end{eqnarray}
    
    Subtracting $f(x)$ gives
    \begin{eqnarray}\label{ineq18:pro1:SecondOrderChar_Dini}
    0\leq\frac{h^2}{2}\max\left\{D^{\prime\prime}f(x;d)+\frac{2o_+(h^2)}{h^2},\,D^{\prime\prime}f(x;-d)+\frac{2o_-(h^2)}{h^2}
    \right\}+\frac{\rho}{2}h^2\|d\|^2.
    \end{eqnarray}
    
    Dividing by $h^2/2$,
    \begin{eqnarray}\label{ineq19:pro1:SecondOrderChar_Dini}
    0\leq\max\left\{D^{\prime\prime}f(x;d)+\frac{2o_+(h^2)}{h^2},\,D^{\prime\prime}f(x;-d)+\frac{2o_-(h^2)}{h^2}\right\}+\rho\|d\|^2.
    \end{eqnarray}
    
    Letting $h\to0^+$ yields
    \begin{eqnarray}\label{ineq20:pro1:SecondOrderChar_Dini}
    \max\left\{D^{\prime\prime}f(x;d),D^{\prime\prime}f(x;-d)\right\}+\rho\|d\|^2\geq 0.
    \end{eqnarray}
\end{proof}

\begin{corollary}\label{cor1:SecondOrderChar_Dini}
 Let $C$ be a nonempty, open and convex subset of $\mathbb{R}^n$ and $f: C \rightarrow\mathbb{R}$ be a $\mathcal{C}^2(C)$ function and $\rho$-weakly quasi-convex function. Then for $x\in C$ and $d\in\mathbb{R}^n$ such that $\left\langle \nabla f(x), d\right\rangle=0$ implies
    \begin{equation}\label{ineq1:cor1:SecondOrderChar_Dini}
        \begin{aligned}
            \left\langle\nabla^2 f(x)d, d\right\rangle+\rho\|d\|^2\geq 0.
        \end{aligned}
    \end{equation}
\end{corollary}

\begin{proof}
    Since $f\in \mathcal{C}^2(C)$, $t\mapsto f(x+td)$ is absolutely continuous on every compact subinterval of $I_{x,d}$. Moreover, 
    \begin{eqnarray}
        f^\prime (x,\pm d)=\pm\left\langle \nabla f(x), d\right\rangle=0,
    \end{eqnarray}
    and
    \begin{eqnarray}
        D^{\prime\prime}f(x;d)=D^{\prime\prime}f(x;-d)=\left\langle\nabla^2 f(x)d, d\right\rangle.
    \end{eqnarray}
    The assertion therefore follows immediately from Proposition \ref{pro1:SecondOrderChar_Dini}.
\end{proof}

\begin{remark}\label{rmk:SecondOrderChar}
    If $f$ is a quasi-convex function then Corollary \ref{cor1:SecondOrderChar_Dini} gives the necessary characterization of it, studied in \cite[Lemma 6.2]{MAVR1972}.
\end{remark}

\begin{corollary}\label{cor2:SecondOrderChar_Dini}
    Let $C$ be a nonempty, open and convex subset of $\mathbb{R}^n$, and let $f:C\rightarrow\mathbb{R}$ be a $\mathcal{C}^2(C)$ function.  Let $f$ be $\rho$-weakly quasi-convex on $C$. Then for $x\in C$ such that $\nabla f(x)=0$ implies
    \begin{equation}\label{ineq1:cor2:SecondOrderChar_Dini}
        \begin{aligned}
            \left\langle\nabla^2 f(x)d, d\right\rangle+\rho\|d\|^2\geq 0,
        \end{aligned}
    \end{equation}
    for all $d\in\mathbb{R}^n.$
\end{corollary}

\begin{proof}
    The proof follows directly from Corollary \ref{cor1:SecondOrderChar_Dini}.
\end{proof}

\begin{corollary}\label{cor3:SecondOrderChar_Dini}
    Let $Q$ be a symmetric matrix of order $n$ and the quadratic function $q(x)=\frac{1}{2}\left\langle Qx,x\right\rangle+\left\langle p,x\right\rangle$, where $p\in\mathbb{R}^n$, is $\rho$-weakly quasi-convex on open convex subset $C$ of $\mathbb{R}^n$. Then for $x\in C$ and $d\in\mathbb{R}^n$ 
   \begin{eqnarray}\label{pro1:QuadFun:ineq1}
       \begin{aligned}
           \left\langle Qx+p, d\right\rangle=0\implies\left\langle Qd,d\right\rangle\geq-\rho\|d\|^2.
       \end{aligned}
   \end{eqnarray}
\end{corollary}

\begin{proof}
    The proof follows directly from Corollary \ref{cor1:SecondOrderChar_Dini} using the facts that for any $x\in C$, $\nabla q(x)=Qx+p$ and $\nabla^2 q(x)=Q.$ 
\end{proof}

\begin{remark}
    When $\rho=0$, Corollary \ref{cor3:SecondOrderChar_Dini} provides a necessary condition for the quadratic function to become quasi-convex on an open subset $C$ of $\mathbb{R}^n.$ 
\end{remark}



The following example shows the converse of Corollary \ref{cor1:SecondOrderChar_Dini} is not true with modulus $\rho$.
\begin{example}\label{example1:cor1:SecondOrderChar_Dini}
    Let $C=\mathbb{R}\times (-\pi/2, \pi/2)\subset\mathbb{R}^2$ and $\rho=1$. Define a function $f:C\rightarrow\mathbb{R}$ by
    \begin{eqnarray}
        f(x)= x_1+\log(\cos x_2),
    \end{eqnarray}
    where $x=(x_1,x_2)$. Then $\nabla f(x)=(1,-\tan x_2)^T$ and 
    \[\nabla^2 f(x) =\begin{bmatrix}
    0 & 0 \\
    0 & -\sec^2 x_2
    \end{bmatrix}.
    \]
    Let $0\neq d=(d_1,d_2)\in\mathbb{R}^2$ such that $\left\langle \nabla f(x), d\right\rangle=0$. This implies $d_1=d_2 \tan x_2$. Now $\left\langle \nabla^2 f(x)d, d\right\rangle+\|d\|^2=d_2^2+d_1^2-d_2^2 \sec^2 x_2=d_2^2 \tan^2 x_2+d_2^2-d_2^2 \sec^2 x_2=0$. Let $x=(1,1)\in C$ and $y=(1,-1)\in C$. Then $f(x)=1+\log(\cos 1)$, $f(y)=1+\log(\cos 1)$, and therefore, $f(x)=f(y)$. Also, $\left\langle \nabla f(x), y-x\right\rangle=2 \tan 1=3.1148$ and $(1/2)\|y-x\|^2=2$. Therefore, $\left\langle \nabla f(x), y-x\right\rangle>(1/2)\|y-x\|^2$. By Proposition \ref{pro1:FirstOrderChar}, $f$ is not $1$-weakly quasi-convex on $C$.
\end{example}

\begin{remark}
\begin{enumerate}
    \item Example \ref{example1:cor1:SecondOrderChar_Dini} demonstrates that a function satisfying the condition in \eqref{ineq1:cor1:SecondOrderChar_Dini} with modulus $\rho=1$, does not guarantee that the function is $1$-weakly quasi-convex on $C$.
    

    \item Since the converse of Proposition \ref{pro1:FirstOrderChar} holds with the weakly quasi-convex modulus $3\rho/2$, it is natural to ask whether an analogous converse of Proposition \ref{pro1:SecondOrderChar_Dini} also holds. That is, let $C$ be a nonempty, open and convex subset of $\mathbb{R}^n$ and $f: C \rightarrow\mathbb{R}$ be a function.
     Assume that 
     \begin{enumerate}
         \item[(A1)] for $x\in C$ and $d\in\mathbb{R}^n$ the directional derivative $f^\prime(x;d)$ exists and $f$ has second-order Dini derivative on $C$.
         \item[(A2)]  for $x\in C$ and $d\in\mathbb{R}^n$ such that $t\mapsto f(x+td)$ is absolutely continuous on every compact interval of 
             \begin{eqnarray*}
                 I_{x,d}=\{t\in \mathbb{R}\;\vert\; x+td\in C\}.
             \end{eqnarray*}
         \item[(A3)] $f^\prime(x,d)=0=f^\prime(x,-d)\implies\max\{D^{\prime\prime}f(x;d),D^{\prime\prime}f(x;-d)\}+\rho\|d\|^2\geq 0$.
     \end{enumerate}
     Does it then follow that $f$ is $3\rho/2$-weakly quasi-convex on $C$?
     The answer to this question is negative. Consider all the data mentioned in Example \ref{example1:cor1:SecondOrderChar_Dini}. Then, for $x=(1,1)\in C$ and $y=(1,-1)\in C$, $f(x)=1+\log(\cos 1)$, $f(y)=1+\log(\cos 1)$, and therefore, $f(x)=f(y)$. Also, $\left\langle \nabla f(x), y-x\right\rangle=2 \tan 1=3.1148$ and $(3/4)\|y-x\|^2=3$. Therefore, $\left\langle \nabla f(x), y-x\right\rangle>(3/4)\|y-x\|^2$. By Proposition \ref{pro1:FirstOrderChar}, $f$ is not $3/2$-weakly quasi-convex on $C$.
\end{enumerate}
\end{remark}


\section{Further properties and examples}
\label{Sec:4}



\begin{proposition}\label{pro2}
    Let $A:\mathbb{R}^n\rightarrow\mathbb{R}^n$ be a linear operator and $\gamma> 0$ be a constant. Let $f:\mathbb{R}^n\rightarrow\mathbb{R}$ be a weakly quasi-convex function with modulus $\rho>0$. The following properties hold.
    \begin{enumerate}
        \item[(i)] $\gamma f$ is weakly quasi-convex with modulus $\gamma\rho$.\label{P1}
        \item[(ii)] $g=f\circ A$ is weakly quasi-convex with modulus $\rho\|A\|^2$.\label{P2}
        \item[(iii)] If $f_i:\mathbb{R}^n\rightarrow\mathbb{R}$ is a weakly quasi-convex function with modulus $\rho_i>0$, then $f=\max_{1\leq i\leq\ell} f_i$ is weakly quasi-convex with modulus $\rho=\max_{1\leq i\leq \ell}\rho_i$.\label{P3}
    \end{enumerate}
\end{proposition}

\begin{proof}
The proofs follow directly from the definition of weakly quasi-convex function.
\end{proof}

\begin{proposition}
    Let $C$ be a nonempty convex subset of $\mathbb{R}^n$, and let $f:C\rightarrow\mathbb{R}$ be a convex function and $L$-Lipschitz continuous.  Assume that the function $\varphi:\mathbb{R}\rightarrow\mathbb{R}$ is $\rho$-weakly quasi-convex and non-decreasing on an interval $I$ containing $f(C)$. Then $\varphi\circ f$ is $\rho L^2$-weakly quasi-convex on $C$. 
\end{proposition}

\begin{proof}
    Let $x,y\in C$. Then $\lambda y+(1-\lambda)x\in C$ for any $\lambda\in [0,1].$

Using our hypotheses, we obtain
    \begin{eqnarray}
        \begin{aligned}
            \varphi\circ f(\lambda y+(1-\lambda)x) &\leq \varphi(\lambda f(y)+(1-\lambda)f(x))\\&
            \leq \max\{\varphi\circ f(y), \varphi\circ f(x)\}+\frac{\rho}{2}\lambda(1-\lambda)\vert f(y)-f(x)\vert^2\\&
            \leq \max\{\varphi\circ f(y), \varphi\circ f(x)\}+\frac{\rho L^2}{2}\lambda(1-\lambda)\|y-x\|^2.
        \end{aligned}
    \end{eqnarray}
\end{proof}

\begin{remark}\label{rmk:pro2}
    In general, the addition or subtraction of two or a finite number of weakly quasi-convex functions is not weakly quasi-convex. 
    A couple of example is provided for the illustration.
\end{remark}

\begin{example}
Define two functions $f_1, f_2:\mathbb{R}\to\mathbb{R}$ by

\begin{eqnarray}
\begin{aligned}
    f_1(x)&=
        \begin{cases}
            0,\;&x<\frac{1}{2}\\
            \frac{1}{2},\; &x\geq \frac{1}{2}
        \end{cases}
    &\text{and}\;\;\;
    f_2(x)&=
        \begin{cases}
            0,\;&x\leq\frac{1}{2}\\
            -\frac{1}{2},\; &x> \frac{1}{2}.
        \end{cases} 
\end{aligned}
\end{eqnarray}

    This implies that
    \begin{eqnarray}
        g(x)=f_1(x)+f_2(x)=
        \begin{cases}
            0,\;&x<\frac{1}{2}\\
            \frac{1}{2},\;&x=\frac{1}{2}\\
            0,\; &x> \frac{1}{2}.
        \end{cases}
    \end{eqnarray}  
    Since $f_1$ and $f_2$ are monotonic on $\mathbb{R}$, they are quasi-convex on $\mathbb{R}$. Therefore, both functions are weakly quasi-convex on $\mathbb{R}$ with any positive modulus.
    Assume that there exists a $\rho>0$ such that for all $x,y\in\mathbb{R}$ and $\lambda\in [0,1]$
    \begin{eqnarray}
       g(\lambda y+(1-\lambda)x)\leq\max\{g(y),g(x)\}+\frac{\rho}{2}\lambda(1-\lambda)\vert y-x\vert^2.
    \end{eqnarray}
    Choose $0<\varepsilon<\sqrt{\frac{2}{\rho}}$.
    By taking 
    \begin{eqnarray}
         x=\frac{1}{2}+\varepsilon,\;\; y=\frac{1}{2}-\varepsilon,\;\;\lambda=\frac{1}{2}.
    \end{eqnarray}
    from the above inequality we get
    \begin{eqnarray}
        \varepsilon\geq\sqrt{\frac{2}{\rho}},
    \end{eqnarray}
    a contradiction.
\end{example}

\begin{example}
Define two functions $f_1,f_2:\mathbb{R}\to\mathbb{R}$ by

\begin{eqnarray}
\begin{aligned}
    f_1(x)&=
        \begin{cases}
            0,\;&x<\frac{1}{2}\\
            \frac{1}{2},\; &x\geq \frac{1}{2}
        \end{cases}
    &\text{and}\;\;\;
    f_2(x)&=
        \begin{cases}
            0,\;&x\leq\frac{1}{2}\\
            \frac{1}{3},\; &x> \frac{1}{2}.
        \end{cases}
\end{aligned}
\end{eqnarray}

    This implies that
    \begin{eqnarray}
        g(x)=f_1(x)-f_2(x)=
        \begin{cases}
            0,\;&x<\frac{1}{2}\\
            \frac{1}{2},\;&x=\frac{1}{2}\\
            \frac{1}{6},\; &x> \frac{1}{2}.
        \end{cases}
    \end{eqnarray}  
    Since $f_1$ and $f_2$ are monotone and quasi-convex on $\mathbb{R}$, both functions are weakly quasi-convex on $\mathbb{R}$ with any positive modulus. Assume that there exists a $\rho>0$ such that for all $x,y\in\mathbb{R}$ and $\lambda\in [0,1]$
    \begin{eqnarray}
        g(\lambda y+(1-\lambda)x)\leq\max\{g(y),g(x)\}+\frac{\rho}{2}\lambda(1-\lambda)\vert y-x\vert^2.
    \end{eqnarray}
    Choose $0<\varepsilon<\sqrt{\frac{2}{3\rho}}$.
    By taking 
    \begin{eqnarray}
         x=\frac{1}{2}+\varepsilon,\;\; y=\frac{1}{2}-\varepsilon,\;\;\lambda=\frac{1}{2},
    \end{eqnarray}
    from the above inequality we get
    \begin{eqnarray}
        \varepsilon\geq\sqrt{\frac{2}{3\rho}},
    \end{eqnarray}
    a contradiction.
\end{example}


Below we provide a sufficient condition ensuring the addition of weakly quasi-convex functions remains weakly quasi-convex.
\begin{lemma}\label{pro2:2}
Let $C$ be a nonempty convex subset of $\mathbb{R}^n$.
    Let $f_1, f_2:C\rightarrow\mathbb{R}$ be two weakly quasi-convex functions on $C$, with modulus $\rho_1, \rho_2>0$, respectively. Assume that for every $x,y\in C$ either $f_1(y)\leq f_1(x)$ and $f_2(y)\leq f_2(x)$ or $f_1(x)\leq f_1(y)$ and $f_2(x)\leq f_2(y)$. Then $g:=f_1+f_2$ is $\rho_1+\rho_2$-weakly quasi-convex on $C$.
\end{lemma}

\begin{proof}
The proof follows directly from the definition of weakly quasi-convex function.
\end{proof}






\begin{lemma}\label{lemma1:conv_sum}
Let $C$ be a nonempty convex subset of $\mathbb{R}^n$.
    Let $f_i:C\rightarrow\mathbb{R}$ be two weakly quasi-convex functions on $C$ for $1\leq i\leq \ell$, with modulus $\rho_i>0$, respectively. Assume that for every  $x, y\in C$ either $f_i(y)\leq f_i(x)$ or $f_i(x)\leq f_i(y)$, for all $1\leq i\leq \ell$. Then $g:=\sum_{i=1}^\ell \lambda_i f_i$ is $\sum_{i=1}^\ell \lambda_i\rho_i$-weakly quasi-convex on $C$, where $\lambda_i\in [0,1]$ such that $\sum_{i=1}^\ell \lambda_i =1$.
\end{lemma}

\begin{proof}
    The proof follows directly from the definition of weakly quasi-convex function.
\end{proof}

\begin{proposition}\label{pro3}
    A function $f:\mathbb{R}^n\rightarrow\mathbb{R}$ is $\rho$-weakly quasi-convex if and only if 
    \begin{eqnarray}
        \begin{aligned}
            g(t)=f\Big(y+\frac{t}{\|x-y\|}(x-y)\Big)
        \end{aligned}
    \end{eqnarray}
    is weakly quasi-convex with modulus $\rho$ on $[0,\|x-y\|]$ for all $x,y\in\mathbb{R}^n$ with $x\neq y.$
\end{proposition}

\begin{proof}
    Let $f$ be $\rho$-weakly quasi-convex on $\mathbb{R}^n$. Fix $t,t^\prime\in[0,\|x-y\|]$ for $x,y\in\mathbb{R}^n$ with $x\neq y.$ Then for every $\lambda\in [0,1]$, we have
    \begin{eqnarray}
        \begin{aligned}
        & g(\lambda t+(1-\lambda)t^\prime) = f\Big(y+\frac{\lambda t+(1-\lambda)t^\prime}{\|x-y\|}(x-y)\Big)\\&
            =f\Big(\lambda\Big[\Big(1-\frac{t}{\|x-y\|}\Big)y+\frac{t}{\|x-y\|}x\Big)\Big]+(1-\lambda)\Big[\Big(1-\frac{t^\prime}{\|x-y\|}\Big)y+\frac{t^\prime}{\|x-y\|}x\Big)\Big]\Big)\\&
            \leq \max\Big\{f\Big(\Big(1-\frac{t}{\|x-y\|}\Big)y+\frac{t}{\|x-y\|}x\Big), f\Big(\Big(1-\frac{t^\prime}{\|x-y\|}\Big)y+\frac{t^\prime}{\|x-y\|}x\Big)\Big\}\\&
            +\frac{\rho}{2}\lambda(1-\lambda)\Big\|t\frac{x-y}{\|x-y\|}-t^\prime\frac{x-y}{\|x-y\|}\Big\|^2\\&
            \leq\max\{g(t),g(t^\prime)\}+\frac{\rho}{2}\lambda(1-\lambda)\vert t-t^\prime\vert^2.
        \end{aligned}
    \end{eqnarray}

  Conversely, suppose $g$ is weakly quasi-convex with modulus $\rho$ on $[0,\|x-y\|]$ for any $x,y\in\mathbb{R}^n$ with $x\neq y.$ Let $s=\frac{t}{\|x-y\|}$, for $0\leq t\leq \|x-y\|$ and for fixed but arbitrary $x,y$.

  \begin{eqnarray}
      \begin{aligned}
          f((1-s) y+ s x)& =g(s\|x-y\|)\\&
          \leq \max\{g(0), g(\|x-y\|)\}+\frac{\rho}{2}s(1-s)\|x-y\|^2\\&
          =\max\{f(y), f(x)\}+\frac{\rho}{2}s(1-s)\|x-y\|^2.
      \end{aligned}
  \end{eqnarray}

  This completes the proof.
\end{proof}


\begin{remark}
    Let $C$ be a open convex subset of $\mathbb{R}^n$. Every $L$-smooth function $f:C\rightarrow\mathbb{R}$  is $L$-weakly convex on C since $f(x)+(L/2)\|x\|^2$ is convex on $C$. Therefore, it is $L$-weakly quasi-convex on $C.$ 
\end{remark}

\subsection{Quadratic functions}

\begin{proposition}\label{pro1:Quadfun}
   Let $Q$ be a symmetric matrix of order $n$ and $\rho>0$ be a constant. Then the quadratic function $q(x)=\frac{1}{2}\left\langle Qx,x\right\rangle+\left\langle p,x\right\rangle$, where $p\in\mathbb{R}^n$, is $\rho$-weakly quasi-convex on $\mathbb{R}^n$ if and only if 
   \begin{eqnarray}\label{pro1::Quadfun:ineq1}
       \begin{aligned}
           \left\langle Qx,x\right\rangle\geq -\rho\|x\|^2,\;\;\text{for all}\;\; x\in\mathbb{R}^n.
       \end{aligned}
   \end{eqnarray}
\end{proposition}
    
\begin{proof} We first prove the sufficiency part of the Proposition. Indeed, consider $z_\lambda=\lambda y+(1-\lambda)x$ for any $\lambda\in (0,1)$ and $x,y\in\mathbb{R}^n$.
    Then we get the following:
    \begin{eqnarray}\label{pro1:QuadFun:ineq2}
        \begin{aligned}
            q(z_\lambda)=\frac{\lambda^2}{2}\left\langle Qy,y\right\rangle &+\frac{(1-\lambda)^2}{2}\left\langle Qx,x\right\rangle+\frac{\lambda(1-\lambda)}{2}\{\left\langle Qx,y\right\rangle+\left\langle Qy, x\right\rangle\}\\&+\lambda\left\langle p,y\right\rangle+(1-\lambda)\left\langle p,x\right\rangle.
        \end{aligned}
    \end{eqnarray}

By our assumptions, we get
      \begin{eqnarray}\label{pro1:QuadFun:ineq3}
        \begin{aligned}
           &\frac{\rho}{2}\lambda(1-\lambda)\|y-x\|^2+\max\Big\{\frac{1}{2}\left\langle Qy,y\right\rangle+ \left\langle p,y\right\rangle,\frac{1}{2}\left\langle Qx,x\right\rangle+ \left\langle p,x\right\rangle\Big\}-q(z_\lambda)\\&\hspace{0.5cm}\geq\frac{\rho}{2}\lambda(1-\lambda)\|y-x\|^2+\frac{\lambda}{2}\left\langle Qy,y\right\rangle+ \lambda\left\langle p,y\right\rangle+\frac{1-\lambda}{2}\left\langle Qx,x\right\rangle\\&
           \hspace{1cm}+ (1-\lambda)\left\langle p,x\right\rangle-\frac{\lambda^2}{2}\left\langle Qy,y\right\rangle-\frac{(1-\lambda)^2}{2}\left\langle Qx,x\right\rangle\\&\hspace{1cm}-\frac{\lambda(1-\lambda)}{2}\{\left\langle Qx,y\right\rangle+\left\langle Qy, x\right\rangle\}-\lambda\left\langle p,y\right\rangle-(1-\lambda)\left\langle p,x\right\rangle\\&
            \hspace{0.5cm}= \frac{\rho}{2}\lambda(1-\lambda)\|y-x\|^2+\frac{1}{4}\lambda(1-\lambda)\left\langle (Q+Q^*)(y-x),y-x\right\rangle\\&\hspace{0.5cm}
            =\frac{1}{2}\lambda(1-\lambda)\{\rho\|y-x\|^2+\left\langle Q(y-x),y-x\right\rangle\}
            \geq 0.
        \end{aligned}
    \end{eqnarray}

    For the necessity, let $x\in\mathbb{R}^n$ be fixed. Denote $\gamma:=\left\langle Qx,x\right\rangle$ and $\delta:=\left\langle p,x\right\rangle$. Without loss of generality, we assume that $x\neq0$ and $\gamma\neq 0$. Consider the function $r(t):=q(t x)$ for $t\in\mathbb{R}$. Therefore, 
    \begin{eqnarray}\label{pro2:QuadFun:ineq2}
        \begin{aligned}
            r(t)=\frac{\gamma}{2}t^2+\delta t.
        \end{aligned}
    \end{eqnarray}
    
    Let $u=(-\delta/\gamma+\varepsilon)x$ and $v=(-\delta/\gamma-\varepsilon)x$, for some constant $\varepsilon\neq0$. Then 
    \begin{eqnarray}\label{pro2:QuadFun:ineq3}
        \begin{aligned}
            & r\left(-\frac{\delta}{\gamma}+\varepsilon\right)-r\left(-\frac{\delta}{\gamma}-\varepsilon\right)\\&=\frac{\gamma}{2}\left(-\frac{\delta}{\gamma}+\varepsilon\right)^2+\delta \left(-\frac{\delta}{\gamma}+\varepsilon\right)-\frac{\gamma}{2}\left(-\frac{\delta}{\gamma}-\varepsilon\right)^2-\delta \left(-\frac{\delta}{\gamma}-\varepsilon\right)\\&
            =0.
        \end{aligned}
    \end{eqnarray}
    
    This implies that $q(u)=r(-\delta/\gamma+\varepsilon)=r(-\delta/\gamma-\varepsilon)=q(v).$ Using Proposition \ref{pro1:FirstOrderChar}, we get
    \begin{eqnarray}\label{pro2:QuadFun:ineq4}
        \begin{aligned}
            \left\langle\nabla q(u), v-u\right\rangle
            \leq\frac{\rho}{2}\|u-v\|^2.
        \end{aligned}
    \end{eqnarray}

    Combining the fact that $\left\langle\nabla q(u), v-u\right\rangle=-2\gamma \varepsilon^2$ and \eqref{pro2:QuadFun:ineq4}, we obtain
    \begin{eqnarray}\label{pro2:QuadFun:ineq6}
        \begin{aligned}
-2 \gamma \varepsilon^2 \leq 2\rho \varepsilon^2 \|x\|^2 ~ \Longleftrightarrow ~ \gamma \geq -\rho\|x\|^2,
        \end{aligned}
    \end{eqnarray}
    that is,
    \begin{eqnarray}\label{pro2:QuadFun:ineq8}
        \begin{aligned}
            \left\langle Qx,x\right\rangle\geq -\rho\|x\|^2.
        \end{aligned}
    \end{eqnarray}
\end{proof}

\begin{corollary}\label{cor1:pro1:QuadFun}
    Let $Q$ be a symmetric matrix of order $n$ and $\rho>0$ be a constant. Let $q(x)=\frac{1}{2}\left\langle Qx,x\right\rangle+\left\langle p,x\right\rangle$, where $p\in\mathbb{R}^n$, be the quadratic function. The following statements are equivalent:
    \begin{enumerate}
        \item[(i)] $q$ is $\rho$-weakly convex on $\mathbb{R}^n$.
        \item[(ii)] $q$ is $\rho$-weakly quasi-convex on $\mathbb{R}^n$.
        \item[(iii)] $Q+\rho I$ is psd.
    \end{enumerate}
\end{corollary}

\begin{remark}\label{pro1:QuadFun:rmk1}
\begin{enumerate}
    \item Note that $q$ is $\rho$-weakly quasi-convex on $\mathbb{R}^n$ if and only if it is $\rho$-weakly convex on $\mathbb{R}^n$.
       However, on a restrictive domain $D$ we may get a quadratic function which is $\rho$-weakly quasi-convex but not $\rho$-weakly convex for some $\rho>0$. For example, define $f(x)=-x^2$ on $\mathbb{R}_{++}$. Since $f$ is monotonically decreasing, it is quasi-convex. Therefore, it is $1$-weakly quasi-convex on $\mathbb{R}_{++}$. However, it is not $1$-weakly convex since $f(x)+(x^2/2)=-x^2/2$ not convex on $\mathbb{R}_{++}$.
    \item In particular, if $\rho=0$, then the above result provides the necessary and sufficient result for the quasi-convexity of the quadratic function (see, for instance, \cite[p. 161]{ABEC2023}). 
\end{enumerate}
\end{remark}

\begin{remark}\label{rmk2:pro1:Quadfun}
$q(x)$ is $\rho$-weakly quasi-convex on an open convex set $C$ if and only  if for any $x,y\in C$ and $\lambda\in (0,1)$
    {\allowdisplaybreaks
    \begin{eqnarray}
    q(\lambda y+(1-\lambda)x)&=&\lambda q(y)+ (1-\lambda) q(x)-\frac{\lambda(1-\lambda)}{2}\left\langle Q(y-x),y-x\right\rangle\nonumber\\
    &\leq& \max\{q(y),q(x)\}+\frac{\rho}{2}\lambda(1-\lambda)\|y-x\|^2.
    \end{eqnarray}
    }

    Since $\lambda\in (0,1)$, the above inequality is true if and only if
    {\allowdisplaybreaks
    \begin{eqnarray}
    q(y)\leq q(x)\implies q(y)-q(x)\leq\frac{1-\lambda}{2}\left\langle (Q+\rho I)(y-x),y-x\right\rangle,
    \end{eqnarray}
    }
    and the above inequality is true if and only if
    {\allowdisplaybreaks
    \begin{eqnarray}
    q(x)\leq q(y)\implies q(x)-q(y)\leq\frac{\lambda}{2}\left\langle (Q+\rho I)(y-x),y-x\right\rangle.
    \end{eqnarray}
    }

    This implies that $q$ is $\rho$-weakly quasi-convex on $C$ if and only if for any $x,y\in C$
    \begin{eqnarray}
       q(x)=q(y)\implies \left\langle (Q+\rho I)(y-x),y-x\right\rangle\geq 0.
    \end{eqnarray}
\end{remark}



\begin{proposition}\label{pro11:Quadfun}
   Let $Q$ be a symmetric matrix of order $n$ and $\rho>0$ be a constant.
   Let a non $\rho$-weakly convex quadratic function $q(x)=\frac{1}{2}\left\langle Qx,x\right\rangle+\left\langle p,x\right\rangle$, where $p\in\mathbb{R}^n$, is $\rho$-weakly quasi-convex on $\mathbb{R}^n_{++}$. Then $Q+\rho I$ has exactly one negative eigenvalue.

\end{proposition}

\begin{proof}
    The quadratic function $q$ is $\rho$-weakly convex on $\mathbb{R}^n_{++}$ if and only if $h(x):=q(x)+(\rho/2)\|x\|^2$ is convex on $\mathbb{R}^n_{++}$. Since $q$ is not $\rho$-weakly convex, $h$ is not convex implies $Q+\rho I$ is not positive semidefinite. Again since $Q+\rho I$ is symmetric, it has at least one negative eigenvalue. 
    
    Since $Q$ is a real symmetric matrix, it has an eigenvalue. Moreover, all eigenvalues are real and there exists a matrix $D$ such that $Q=PDP^T$ by the spectral theorem, where $D$ is a diagonal matrix with diagonal elements that are the eigenvalues of $Q$ and $P$ is the matrix with orthonormal eigenvectors corresponding to the eigenvalues of $Q$.
    
    Let $\{\lambda_1,\lambda_2,\dots, \lambda_n\}$ be the set of eigenvalues and $\{v_1,v_2,\dots,v_n\}$ be the set of orthonormal eigenvectors corresponding to each eigenvalue of $Q$. Therefore, $\{\lambda_1+\rho ,\lambda_2+\rho,\dots, \lambda_n+\rho\}$ is the set of eigenvalues and $\{v_1,v_2,\dots,v_n\}$ is the set of orthonormal eigenvectors corresponding to each eigenvalue of $Q+\rho I$. 

    Suppose to the contrary that it has more than one, without loss of generality, we may assume that it has at least two negative eigenvalues (counted with multiplicity), say $\lambda_1+\rho$ and $\lambda_2+\rho$. Consider $d\in\text{span}\{v_1,v_2\}$ of the form $d= r v_1+t v_2$, where $r, t$ are scalars. Let $x\in\mathbb{R}^n_{++}$, $\alpha=\left\langle Qx+p,v_1\right\rangle\neq 0$, and $\beta=\left\langle Qx+p,v_2\right\rangle\neq 0$. Then 
    \begin{eqnarray}
        \left\langle Qx+p,d\right\rangle
        = r\left\langle Qx+p,v_1\right\rangle+t\left\langle Qx+p,v_2\right\rangle=0,
    \end{eqnarray}   
    implies that $r\alpha+t\beta=0$, the linear equation has a nonzero solution. For example, take $t=-\alpha$ and $r=\beta$. 
    When $(\alpha,\beta)\neq (0,0)$, 
    \begin{eqnarray}
        \left\langle Qd,d\right\rangle+\rho\|d\|^2
        &=& \left\langle r \lambda_1 v_1+t \lambda_2 v_2,r v_1+t v_2\right\rangle+\rho\left\langle r v_1+t v_2,r v_1+t v_2\right\rangle\\
        &=& r^2(\lambda_1+\rho)+t^2(\lambda_2+\rho)<0,\nonumber
    \end{eqnarray}     
    a contradiction to Corollary \ref{cor3:SecondOrderChar_Dini}. 
    Similarly, when $(\alpha,\beta)=(0,0)$, choose $d\in\text{span}\{v_1,v_2\}$. 
    In particular, $d=v_1$ implies that
    \begin{eqnarray}
        \left\langle Qd,d\right\rangle+\rho\|d\|^2=\lambda_1+\rho<0,\nonumber
    \end{eqnarray}    
    a contradiction to Corollary \ref{cor3:SecondOrderChar_Dini}.
    \end{proof}

\begin{remark}\label{rmk1:pro11:Quadfun}
The converse of Proposition \ref{pro11:Quadfun} is not true. For example, take $n=2$, $\rho=1$, and
\[Q=\begin{bmatrix}
    -2 & 0 \\
    0 & 0
\end{bmatrix}
\] with $p=(0, 1)$. Then $q(x)=-x_1^2+x_2$ on $\mathbb{R}^2_{++}$ and $Q+I$ has exactly one negative eigenvalue. Let $x=(1/4,1)$ and $d=(1,1/2)$. Then $\left\langle Qx+p,d\right\rangle= 0$ but $\left\langle (Q+I)d, d\right\rangle=-3/4<0$. Therefore, $q$ is not $1$-weakly quasi-convex on $\mathbb{R}^2_{++}$ by Corollary \ref{cor3:SecondOrderChar_Dini} with $C=\mathbb{R}^2_{++}$.
\end{remark}

\begin{remark}
\begin{enumerate}
    \item[(1)]  A result stated in \cite[Theorem 3]{MVJO1993} for a nonconvex function to become quasi-convex on $\mathbb{R}^n_+$. A similar result does not follow for the class of weakly quasi-convex functions on $\mathbb{R}^n_{++}$. Proposition \ref{pro11:Quadfun} shows that a non $\rho$-weakly convex function $q$ is $\rho$-weakly quasi-convex on $\mathbb{R}^n_{++}$ implies only $Q+\rho I$ has exactly one negative eigenvalue. Therefore, for example, let $n=2$, $\rho=1$, and
    \[Q=\begin{bmatrix}
       -2 & 0\\
       0 & 0
    \end{bmatrix}
    \] with $p=(0, 0)$. Then $q(x)=-x_1^2$ is not $1$-weakly convex on $\mathbb{R}^2_{++}$. 
    However, it is quasi-convex on $\mathbb{R}^2_{++}$. 
    Therefore, in particular, the function $q$ is $1$-weakly quasi-convex on $\mathbb{R}^2_{++}$. Now,
    \[Q+I=\begin{bmatrix}
       -1 & 0\\
       0 & 1
    \end{bmatrix}\] has exactly one negative eigenvalue. However, $Q+I\nless 0$.

    \item[(2)] If $Q+\rho I<0$ then $q$ may not be $\rho$-weakly quasi-convex on $\mathbb{R}^n_{++}$. For example, let 
    $n=2$, $\rho=1$, and
    \[Q=\begin{bmatrix}
       -2 & -1\\
       -1 & -3
    \end{bmatrix}
    \] with $p=(0, 0)$. Then $q(x)=-(1/2)(2x_1^2+2x_1 x_2+3x_2^2)$ on $\mathbb{R}^2_{++}$. Let $x=(1,1)$ and $d=(1,-3/4)$. Then $\left\langle Qx+p,d\right\rangle= 0$ but $\left\langle (Q+I)d, d\right\rangle<0$. Therefore, $q$ is not $1$-weakly quasi-convex on $\mathbb{R}^2_{++}$ by Corollary \ref{cor3:SecondOrderChar_Dini} with $C=\mathbb{R}^2_{++}$.

    Similarly, we can find a non $\rho$-weakly convex function which is $\rho$-weakly quasi-convex on $\mathbb{R}^n_{++}$ but $Q+\rho I<0$ does not hold. For example, take the data presented in $(1)$.
\end{enumerate}

\end{remark}

\begin{remark}
    Let $0\in C$ be a convex cone of $\mathbb{R}^n$. 
    Consider the quadratic function $q(x)=\frac{1}{2}\left\langle Qx,x\right\rangle+\left\langle p,x\right\rangle$, where $p\in\mathbb{R}^n$ and $Q$ is symmetric (non-diagonal) of order $n$. 
    Then also $\rho$-weak quasi-convexity of $q$ does not imply $\rho$-weak convexity. For example take $C=\mathbb{R}^n_+$, $n=2$, $p=(0,0)$, and
    \begin{eqnarray}
        Q=\begin{bmatrix}
            -1 & -2 \\
            -2 & -1
        \end{bmatrix}.
    \end{eqnarray}
    Then $q$ in quasi-convex on $C$, and therefore, it it $\rho$-weakly quasi-convex on $C$ for any $\rho>0$. Let $\rho=1$ and $d=(1,1)$. Define $g(x):=q(x)+(1/2)\|x\|^2$ for $x\in C$. This implies $\left\langle \nabla^2g(x)d,d\right\rangle=-4<0.$ Thus, $g$ is not convex on $C$.
\end{remark}

\section{Optimization consequences}\label{Sec:OptimizationProblem}

\subsection{First-order stationarity and variational inequality}

\begin{proposition}\label{pro1:OptimizationProblem}
    Let $C$ be an open convex subset of $\mathbb{R}^n$. Let $f:C\rightarrow\mathbb{R}$ be $\rho$-weakly quasi-convex and differentiable. Let $\mathcal{Z}$ be a nonempty convex subset of $\mathbb{R}^n$ such that $\mathcal{Z}\subset C$. Let $x^*\in\mathcal{Z}$ be a minimizer of the \textit{constrained optimization problem}
\begin{eqnarray}\label{OptimizationProblem}
    \begin{aligned}
        \min_{x\in\mathcal{Z}} f(x).
    \end{aligned}
\end{eqnarray}
Then 
\begin{eqnarray}\label{ineq1:OptimizationProblem}
        \left\langle\nabla f(x),x-x^*\right\rangle+\frac{\rho}{2}\|x-x^*\|^2\geq 0,\;\;\text{for all}\;\; x\in\mathcal{Z}.
\end{eqnarray}
\end{proposition}

\begin{proof}
    Let $f$ be $\rho$-weakly quasi-convex and differentiable on $C$. By the definition of a minimizer, $f(x^*)\leq f(x)$ for all $x\in \mathcal{Z}$. Then by Proposition \ref{pro1:FirstOrderChar}, we get
    \begin{eqnarray}\label{ineq2:OptimizationProblem}
         \left\langle\nabla f(x),x-x^*\right\rangle+\frac{\rho}{2}\|x-x^*\|^2\geq 0,\;\;\text{for all}\;\; x\in\mathcal{Z}.
    \end{eqnarray}
\end{proof}

\begin{remark}\label{rmk2:OptimizationProblem}
 The inequality \eqref{ineq2:OptimizationProblem} is only a necessary consequence. However, it is not a sufficient condition. Consider $f(x)=-x^3$ for $x\in C=\mathbb{R}$. Let $\mathcal{Z}=(-\infty,1/3]$. The function $f$ is $\rho$-weakly quasi-convex on $C$ for any $\rho>0$. Let $\rho=2$ and $x^*=0$. Then $\left\langle\nabla f(x),x-x^*\right\rangle+\frac{\rho}{2}\|x-x^*\|^2\geq 0$ for any $x\in\mathcal{Z}$ but $0$ is not a minimizer of $f$.
\end{remark}

\begin{remark}\label{rmk3:OptimizationProblem}
\begin{enumerate}
    \item[(1)] In particular, if we set $\rho=0$ in \eqref{ineq1:OptimizationProblem}, then it reduces to the Minty variational inequality associated to $\nabla f$ on $\mathcal{Z}$ (see, for instance, \cite{GJMI1967,BBUL2022}), i.e., find a point $x^*\in \mathcal{Z}$ such that
    \begin{eqnarray}\label{minty:vip}
        \left\langle \nabla f(x),x-x^*\right\rangle\geq 0,\;\;\text{for all}\;\; x\in \mathcal{Z}.
    \end{eqnarray}
    The inequality in \eqref{ineq1:OptimizationProblem}, we call it a \textit{relaxed Minty variational inequality} associated to $\nabla f$ on $\mathcal{Z}$.

    \item[(2)] Assume that 
    \begin{eqnarray}\label{ineq3:OptimizationProblem}
        \left\langle\nabla f(x)-\nabla f(y),x-y\right\rangle\leq -(\rho/2)\|x-y\|^2
    \end{eqnarray}
    for any $x,y\in \mathcal{Z}$ and for some $\rho>0$. That is, if $-\nabla f$ is $\rho/2$-strongly monotone on $\mathcal{Z}$, then \eqref{ineq1:OptimizationProblem} implies $\left\langle\nabla f(x^*),x-x^*\right\rangle\geq 0$ for all $x\in \mathcal{Z}$. This is called the variational inequality problem associated to $\nabla f$ on $\mathcal{Z}$ (see, for instance, \cite{GSTA1964}).

    \item[(3)] Any convex and differentiable function does not satisfy the above inequality $\eqref{ineq3:OptimizationProblem}$. However, a quasi-convex function may satisfy this. For example, take a function $f(x)=-\|x\|^2$ on $C=\mathbb{R}^n_+$ with $\mathcal{Z}=[0,1]^n$, unit cube in $\mathbb{R}^n$. It is clear that $f$ is monotonically decreasing, and hence it is quasi-convex. The function $f$ satisfies the inequality \eqref{ineq3:OptimizationProblem} with $\rho=4$.
\end{enumerate}
\end{remark}

\subsubsection{Application to solving saddle point problems}\label{Sec:saddleProb}
Let $C$ and $Q$ be two nonempty, open and convex subsets of $\mathbb{R}^n$ and $\mathbb{R}^m$, respectively. 
Let $\widetilde{C}\subset C$ and $\widetilde{Q}\subset Q$ be two convex subsets of $\mathbb{R}^n$ and $\mathbb{R}^m$, respectively. 
Let $\psi:C\times Q\to\mathbb{R}$ be a function.
Then the \textit{saddle point problem} consists of finding a point $(x^*,y^*)\in \widetilde{C}\times \widetilde{Q}$ such that
\begin{eqnarray}\label{ineq1:saddleProb}    \min_{x\in\widetilde{C}}\max_{y\in\widetilde{Q}}\psi(x,y)=\psi(x^*,y^*)=\max_{y\in\widetilde{Q}}\min_{x\in\widetilde{C}}\psi(x,y).
\end{eqnarray}
For more details on saddle point problems in different contexts, we refer \cite{MSIO1958,JVON1928}.

\begin{proposition}\label{pro1:saddleProb}
Let $C$ and $Q$ be two nonempty, open and convex subsets of $\mathbb{R}^n$ and $\mathbb{R}^m$, respectively. 
Let $\widetilde{C}\subset C$ and $\widetilde{Q}\subset Q$ be two convex subsets of $\mathbb{R}^n$ and $\mathbb{R}^m$, respectively. 
Let $\psi:C\times Q\to\mathbb{R}$ be a function such that $\psi(\cdot, y)$ is $\rho$-weakly quasi-convex for every $y\in \title{Q}$. Then
\begin{eqnarray}\label{ineq1:pro1:saddleProb}
    \widetilde{\psi}(x)=\max_{y\in\widetilde{Q}}\psi(x,y)
\end{eqnarray}
is weakly quasi-convex on $C$ with modulus $\rho$.
\end{proposition}

\begin{proof}
    Fix $y\in\widetilde{Q}$. Since $\psi(\cdot, y)$ is $\rho$-weakly quasi-convex, for every $x,\omega\in  C$ and every $\lambda\in [0,1]$
    \begin{eqnarray}\label{ineq2:pro1:saddleProb}
        \psi(\lambda\omega+(1-\lambda)x,y)\leq \max\{\psi(\omega,y),\psi(x,y)\}+\frac{\rho}{2}\lambda(1-\lambda)\|x-\omega\|^2.
    \end{eqnarray}
    This implies that
    \begin{eqnarray}\label{ineq3:pro1:saddleProb}
        \max_{y\in\widetilde{Q}}\psi(\lambda\omega+(1-\lambda)x,y)\leq \max\left\{\max_{y\in\widetilde{Q}}\psi(\omega,y),\max_{y\in\widetilde{Q}}\psi(x,y)\right\}+\frac{\rho}{2}\lambda(1-\lambda)\|x-\omega\|^2.
    \end{eqnarray}
    Thus, 
    \begin{eqnarray}\label{ineq4:pro1:saddleProb}
        \widetilde{\psi}(\lambda\omega+(1-\lambda)x)\leq \max\left\{\widetilde{\psi}(\omega),\widetilde{\psi}(x)\right\}+\frac{\rho}{2}\lambda(1-\lambda)\|x-\omega\|^2.
    \end{eqnarray}    
\end{proof}

\begin{proposition}\label{pro2:saddleProb}
Let $C$ and $Q$ be two nonempty, open and convex subsets of $\mathbb{R}^n$ and $\mathbb{R}^m$, respectively. 
Let $\widetilde{C}\subset C$ and $\widetilde{Q}\subset Q$ be two convex subsets of $\mathbb{R}^n$ and $\mathbb{R}^m$, respectively. 
Let $\psi:C\times Q\to\mathbb{R}$ be a function such that $\psi(\cdot, y)$ is $\rho$-weakly quasi-convex, and let
\begin{eqnarray}
    \widetilde{\psi}(x)=\max_{y\in\widetilde{Q}}\psi(x,y)
\end{eqnarray}
be differentiable on $C$. Assume that the saddle point problem \eqref{ineq1:saddleProb} has a solution. Then 
\begin{eqnarray}\label{ineq1:pro2::saddleProb}
        \left\langle\nabla \widetilde{\psi}(x),x-x^*\right\rangle+\frac{\rho}{2}\|x-x^*\|^2\geq 0,\;\text{for all}\; x\in\widetilde{C},
\end{eqnarray}

\end{proposition}

\begin{proof}
    By Proposition \ref{pro1:saddleProb}, we conclude that $\widetilde{\psi}(x)$ is $\rho$-weakly quasi-convex on $C$.    
    Let $(x^*,y^*)$ be a saddle point of the problem \eqref{ineq1:saddleProb} by our hypothesis. 
    Therefore,
    \begin{eqnarray}
    \psi(x^*,y^*)\leq \psi(x,y^*)\leq\max_{y\in\widetilde{Q}}\psi(x,y):=\widetilde{\psi}(x)
    \end{eqnarray}
    for all $x\in\widetilde{C}$. Consequently, $x^*\in\argmin_{x\in\widetilde{C}}\widetilde{\psi}(x)$.
    Using the differentiability hypothesis, the result directly follows from Proposition \ref{pro1:OptimizationProblem}.
\end{proof}

\begin{remark}
    A sufficient condition for the differentiability of $\widetilde{\psi}$ can be obtained from Danskin's theorem \cite[Proposition B.25]{DPBE1999}.
\end{remark}

\subsection{Weak quasi-convexity and approximate solution sets of a constrained minimization problem}

Let $C$ be a nonempty convex subset of $\mathbb{R}^n$.
Let $f:C\to\mathbb{R}$ be a function. Assume that 
\begin{eqnarray}\label{min:2}
    f^*=\inf_{x\in C} f(x)>-\infty.
\end{eqnarray}
That is, the function $f$ is bounded below on $C$.
For any $\varepsilon>0$, the approximate solution set $S_\varepsilon$ is defined as
\begin{eqnarray}
    S_\varepsilon=\{x\in C\;\vert\; f(x)\leq f^*+\varepsilon\},
\end{eqnarray}
which is nonempty by the definition of the infimum. However, $S_0=\argmin_{x\in C}$ f(x) could be empty in general.

\begin{theorem}\label{thm1:approx:soul}
Let $C$ be a nonempty convex subset of $\mathbb{R}^n$.
Let $f:C\to\mathbb{R}$ be a function that is bounded below on $C$. Let $\varepsilon>0$ be arbitrary. Then $f$ is $\rho$-weakly quasi-convex on $C$ if and only if 
\begin{eqnarray}\label{ineq1:thm1:approx:soul}
     x,y\in S_\varepsilon\implies\lambda y+(1-\lambda)x\in S_{\varepsilon+\frac{\rho}{2}\lambda(1-\lambda)\|y-x\|^2},\;\;\text{for all}\;\;\lambda\in [0,1].
\end{eqnarray}
\end{theorem}

\begin{proof}
    Let $x,y\in S_\varepsilon$. Therefore, $f(x)\leq f^*+\varepsilon$ and $f(y)\leq f^*+\varepsilon$. Suppose $f$ is $\rho$-weakly quasi-convex on $C$. Then 
    \begin{eqnarray}\label{ineq2:thm1:approx:soul}
        f(\lambda y+(1-\lambda)x)\leq f^*+\varepsilon+\frac{\rho}{2}\lambda(1-\lambda)\|y-x\|^2.
    \end{eqnarray}
    Conversely, assume that the implication in \eqref{ineq1:thm1:approx:soul} holds. Let $x,y\in C$. Consider 
    \begin{eqnarray}\label{ineq3:thm1:approx:soul}
        \varepsilon_0=\max\{f(y),f(x)\}-f^*\;\;\text{and}\;\;\varepsilon=\varepsilon_0+\delta,
    \end{eqnarray}
    for some $\delta>0$ arbitrarily small.
    Then $f(x)\leq f^*+\varepsilon_0$ and $f(y)\leq f^*+\varepsilon_0$. Thus  $\varepsilon> 0$ and $x,y\in S_\varepsilon$. Using our hypothesis in \eqref{ineq1:thm1:approx:soul}, we obtain
    \begin{eqnarray}\label{ineq4:thm1:approx:soul}
        f(\lambda y+(1-\lambda)x)\leq f^*+\varepsilon+\frac{\rho}{2}\lambda(1-\lambda)\|y-x\|^2,
    \end{eqnarray}
    which further implies by taking $\delta\downarrow 0$
    \begin{eqnarray}\label{ineq5:thm1:approx:soul}
        f(\lambda y+(1-\lambda)x)\leq \max\{f(y),f(x)\}+\frac{\rho}{2}\lambda(1-\lambda)\|y-x\|^2.
    \end{eqnarray}    
\end{proof}

If we set $\rho=0$ in the above Theorem \ref{thm1:approx:soul} we obtain the following corollary concerning quasi-convex functions.
\begin{corollary}\label{cor1:thm1:approx:soul}
Let $C$ be a nonempty convex subset of $\mathbb{R}^n$.
Let $f:C\to\mathbb{R}$ be a function and it is bounded below on $C$. Let $\varepsilon>0$ be arbitrary. Then $f$ is quasi-convex on $C$ if and only if 
\begin{eqnarray}\label{ineq1:cor1:thm1:approx:soul}
     x,y\in S_\varepsilon\implies\lambda y+(1-\lambda)x\in S_{\varepsilon},\;\text{for all }\lambda\in [0,1].
\end{eqnarray}
\end{corollary}

\begin{remark}
    From Corollary \ref{cor1:thm1:approx:soul} it follows that for the class of quasi-convex functions the associated approximate solution set is convex.
\end{remark}

\begin{theorem}\label{thm2:approx:soul}
Let $C$ be a nonempty, open and convex subset of $\mathbb{R}^n$.
Let $f:C\to\mathbb{R}$ be a function and it is bounded below on $C$. Let $\varepsilon>0$ be arbitrary. Suppose that $f$ is $\rho$-weakly quasi-convex and differentiable on $C$. Assume that the boundary of $S_\varepsilon$ is $\partial S_\varepsilon$ such that $\partial S_\varepsilon\cap C\neq \varnothing$. Then for all $x\in S_\varepsilon$ and $x^*\in \partial S_\varepsilon\cap C$
\begin{eqnarray}\label{ineq1:thm2:approx:soul}
\left\langle \nabla f(x^*), x-x^*\right\rangle\leq \frac{\rho}{2}\|x-x^*\|^2.
\end{eqnarray}
\end{theorem}

\begin{proof}
    Let $x\in S_\varepsilon$ and $x^*\in \partial S_\varepsilon\cap C$. Then $f(x)\leq f(x^*)$. The result directly follows from Proposition \ref{pro1:FirstOrderChar}.
\end{proof}

\section{Conclusion}
\label{Sec:6}
The notion of weak quasi-convexity opens a new avenue for extending convex analysis beyond the classical convex and weakly convex settings. 
The results presented in this paper show that the class retains sufficient structure to support a meaningful analytical treatment. 
Additionally, the analysis carried out in this paper for the class of quadratic functions indicates that properties valid on the whole space may change substantially when one passes to restricted convex sets or cones. 
The notions of monotonicity associated with weakly quasi-convex functions generalize various types of existing monotonicity concepts. 
Therefore, the class of weakly quasi-convex functions genuinely provides distinct structure rather than being only a minor relaxation of existing notions. 
Moreover, the open question posed in this paper provides a direction for the theoretical investigation, while the results presented in this paper open a new avenue for the algorithmic study for this class of functions.  We aim to explore additional consequences and applications of weakly quasi-convex functions in optimization. 

\bigskip{}
\textbf{Acknowledgment}
This work was supported in part by the Israeli Science Foundation (ISF), Grant No. 992/23. We would like to thank Professor Felipe Lara, Institute for Advanced Studies, University of Tarapacá, Arica, Chile, for his valuable comments and suggestions, which helped improve an earlier version.





\end{document}